\documentclass[11pt, oneside]{article}   	
\usepackage[margin=1in]{geometry}                		
\usepackage{graphicx}				
\usepackage{amsmath}
\usepackage{amsthm}
\usepackage{mathtools}
\usepackage[dvipsnames]{color}								
\usepackage{amssymb}
\usepackage[round]{natbib}
\usepackage{hyperref}
\usepackage{xcolor}	
\hypersetup{
     colorlinks   = true,
     citecolor     = teal, 
     linkcolor     = teal
}

\usepackage{soul}
\usepackage{endnotes}

\newcommand{\bs}{\boldsymbol}

\newcommand{\transpose}{^\text{T}}
\newcommand{\ddt}{\frac{d}{dt}}

\begin{document}

\begin{center}\Large
	Inaccuracy of Ensemble-Based Covariance Propagation,\\ Beyond Sampling Error

\vspace{0.5em}

\normalsize
Shay Gilpin$^1$\\ \footnotesize \flushleft
$^1$Department of Mathematics, University of Arizona, Tucson, Arizona, USA, sgilpin@arizona.edu\\
\end{center}


\section*{Abstract}
Modern data assimilation schemes typically use the same discrete dynamical model to evolve the state estimate in time also to approximate the evolution, or propagation, of the estimation error covariance. Ensemble-based methods, such as the ensemble Kalman filter, approximate the evolution of the covariance through the propagation of individual ensemble members. Thus, it is tacitly assumed that if the discrete state propagation and resulting mean state estimates are accurate, then the ensemble-based discrete covariance propagation will be accurate as well, apart from sampling errors due to limited ensemble size. Through a series of numerical experiments supported by analytical results, we demonstrate that this assumption is false when correlation length scales approach grid resolution. We show for states that satisfy advective dynamics, that while the discrete state propagation and ensemble mean state estimates are accurate, the corresponding ensemble covariances can be remarkably inaccurate, well beyond that expected from sampling errors or typical numerical discretization errors. The underlying problem is a fundamental discrepancy between discrete covariance propagation and the continuum covariance dynamics, which we can identify because the exact continuum covariance dynamics are known. Errors in the ensemble covariances, which can be at least one order of magnitude larger than those of the mean state when correlation lengths begin to approach grid scale, cannot be rectified by the usual methods, such as covariance inflation and localization. This work brings to light a fundamental problem for data assimilation schemes that propagate covariances using the same discrete dynamical model used to propagate the state.\\

\noindent\textit{Keywords}: covariance propagation, ensemble data assimilation, advective dynamics, Kalman filter

\section{Introduction}
Data assimilation is widely practiced in the atmospheric and ocean sciences to improve estimates of the state of a given dynamical system, and possibly its parameters, by using observation information and characterizing uncertainties. Uncertainty is quantified through the estimation error covariance, and covariance propagation, i.e., the evolution of the estimation error covariance in time, is an important component of modern data assimilation algorithms. 

Generally, covariance propagation is implied by the propagation of the state. In fact, the same discrete dynamical model that is used to propagate the state estimate, or possibly a linearized version of it, is typically used to also approximate the propagation of the covariance matrix. The linearized discrete dynamical model is often represented by a matrix $\mathbf{M}$. Kalman filters, for example, propagate the covariance matrix explicitly using the same state transition matrix $\mathbf{M}$ defined for the state propagation \citep{kalman1960new,jazwinski1970stochastic}. Variational-based schemes that evolve covariances, such as 4D-Var, evolve the covariance matrix implicitly through a tangent linear model $\mathbf{M}$ and its adjoint during the minimization of the cost function \citep{lorenc2003modelling}. Ensemble-based methods, such as the ensemble Kalman filter (EnKF), approximate the evolution of the covariance matrix through the nonlinear evolution of the individual ensemble members \citep{evensen1994sequential,evensen2009ensemble}.

Recent work in chemical constituent data assimilation has observed that full-rank covariance propagation can be quite inaccurate when using the same discrete dynamical model $\mathbf{M}$ defined for discrete state propagation \citep{menard2000aassimilation,menard2021numerical,menard2000bassimilation,lyster2004lagrangian, pannekoucke2021methodology,gilpin2022continuum,gilpin2025inaccuracy}. Propagated covariances suffer from spurious loss and gain of variance and from inaccurate correlations both near and away from zero separation. These inaccuracies in the covariances become particularly pronounced in regions where correlation length scales approach grid resolution \citep{menard2000aassimilation,lyster2004lagrangian,gilpin2022continuum,gilpin2025inaccuracy}. This behavior observed during full-rank covariance propagation raises the question of whether ensemble-based covariance propagation might also suffer from the same inaccuracies. Since the introduction of the EnKF by \cite{evensen1994sequential}, the main focus of effort has traditionally been on understanding and mitigating sampling error, which is due to covariance matrices being estimated from a limited ensemble \citep[e.g.,][]{houtekamer1998data,hamill2001distance,evensen2009ensemble,whitaker2012evaluating}. The covariance propagation itself, via the evolution of the ensemble members, as a source of error has largely been overlooked in the literature.

For states that satisfy the linear advection equation and related hyperbolic partial differential equations (PDEs), we find in this work that ensemble-based covariance propagation is severely inaccurate, particularly when correlation lengths become small. We demonstrate this through a series of numerical experiments and compare these results with previous work that has considered the full-rank, explicit covariance propagation of the Kalman filter \citep{gilpin2022continuum,gilpin2025inaccuracy}. We are able to explicitly quantify the errors in covariance propagation because the exact continuum covariance dynamics are known. The continuum covariance dynamics is governed by a deterministic PDE that can be explicitly formulated and solved in special cases, and has been studied extensively by \cite{cohn1993dynamics} and \cite{gilpin2022continuum}. 

In these numerical experiments, we find that the propagation of the mean state by the ensemble is accurate, while at the same time, the ensemble-based covariance propagation is remarkably inaccurate, much more so than expected from either sampling errors or the usual numerical discretization errors. The errors in the ensemble covariances, which can be at least one order of magnitude larger than errors in the ensemble mean, are relatively constant with respect to ensemble size, but increase rapidly as correlation lengths decrease. The ensemble variances themselves hardly resemble the exact variance, exhibiting  both spurious loss and gain of variance in some cases, or severe variance loss in others, neither of which can be addressed by the usual techniques of variance inflation. The ensemble correlations, while exhibiting some spurious long-range correlations as ensemble sizes decrease, can in fact be severely inaccurate near the covariance diagonal. Methods such as covariance localization \citep{houtekamer2001sequential,hamill2001distance,hamill2009comments} cannot help to rectify this issue. In fact, we find that ensemble covariances approximate covariances propagated explicitly at full rank as in a Kalman filter rather well, but neither can approximate the exact, known covariance dynamics with acceptable accuracy. 

The primary source of the inaccuracy in both ensemble-based and explicit, full-rank covariance propagation, is a fundamental discrepancy between the discrete dynamical model $\mathbf{M}$ that is built from the continuum state dynamics and the behavior of the continuum covariance dynamics. While the state and covariance are related by definition, as shown in \cite{cohn1993dynamics} and \cite{gilpin2022continuum}, the solutions to the continuum covariance PDE behave in ways that are not expected when considering the state dynamics alone. As a consequence, the discrete dynamical model $\mathbf{M}$, which encodes the state dynamics, is unable to capture the correct covariance dynamics when correlation lengths become small. This is observed in our numerical experiments, and is also supported by the analysis of the variance dynamics in \cite{gilpin2025inaccuracy} and the correlation dynamics here in Appen~\ref{sec:appendix correlation}.

Errors in covariance propagation can sometimes be attributed to the notion of ``model error," often represented by a stochastic forcing term that is added to the state propagation, and consequently a model error covariance matrix added to the covariance propagation (\citealp[Ch.~13.3]{daley1991atmospheric}; \citealp{dee1995line}). In the present work, we assume the true dynamical system is given by an exact PDE for the state and covariance so that we can clearly isolate the errors caused by discrete propagation and ensemble size.

Our intention with this work is to clearly demonstrate through a series of simple numerical experiments a fundamental problem with covariance propagation in data assimilation for advective dynamics: accurate state propagation does not imply accurate covariance propagation when using the same discrete dynamical model for both. Ensemble-based covariance propagation is popular because covariance propagation is computationally tractable in these algorithms. As seen in our experiments, however, the impact of using the same discrete dynamical model on the covariances themselves has been overlooked in the literature and can be quite problematic. The goal of this work is thus to clearly present the problem and open the door to potential solutions.

This paper is organized as follows. We begin in Sec.~\ref{sec:background} by reviewing the propagation of the discretized mean state and covariance in data assimilation algorithms, specifically for full-rank propagation of the Kalman filter and its ensemble approximations, such as for the EnKF.  This is followed by the continuum formulation and a series of numerical experiments in Sec.~\ref{sec:experiments}. In these experiments, we compare the mean state and covariance estimated from the propagated ensemble with the exact mean state and covariance to illustrate the severity of the errors incurred during discrete covariance propagation. We then discuss these results in Sec.~\ref{sec:discussion} in the context of the continuum dynamics and analysis presented in \cite{cohn1993dynamics} and \cite{gilpin2022continuum,gilpin2025inaccuracy}, followed by concluding remarks in Sec.~\ref{sec:conclusions}. Appendix \ref{sec:appendix correlation} provides our analysis of the correlation dynamics approximated by the numerical schemes and identifies the error terms that produce the behavior observed in the numerical experiments.

\section{Discrete State and Covariance Propagation}\label{sec:background}
A data assimilation cycle consists of two steps. The first is a time propagation step (forecast) which evolves the system (e.g., the mean state, covariance, ensemble of states) forward in time according to given dynamical system. The second is the measurement update step (analysis) that incorporates the observation information and uncertainty using Bayesian inference to produce an updated state estimate and corresponding uncertainty. The focus of this work is on the time propagation step alone.

Let $\mathbf{q}_k$ represent the $N$-vector discretized state of a dynamical system at time $t_k$. This state is evolved forward to time $t_{k+1}$ according to
\begin{equation}\label{eq:state prop}
	\mathbf{q}_{k+1} = \mathbf{M}_{k,k+1}\mathbf{q}_k,
\end{equation}  
where $ \mathbf{M}_{k,k+1}$ is an $N\times N$ propagation (state transition) matrix that represents the discretized model dynamics. For simplicity, we assume the state dynamics are deterministic and linear, with no forcing, random or otherwise. 

If the initial state $\mathbf{q}_0$ is random with mean $\overline{\mathbf{q}}_0$, then the mean state at time $t_k$, $\overline{\mathbf{q}}_k$, satisfies the same discrete propagation as the state, \eqref{eq:state prop}, since $\mathbf{M}_{k,k+1}$ is deterministic,
\begin{equation}\label{eq:mean state prop}
	\overline{\mathbf{q}}_{k+1} = \mathbf{M}_{k,k+1}\overline{\mathbf{q}}_k.
\end{equation} 
Along with the mean state, we can define the $N\times N$ covariance matrix $\mathbf{P}_k$ at time $t_k$, where the $i,j$-th entry of the covariance matrix is defined as
\begin{equation}\label{eq:covariance matrix def}
	\big(\mathbf{P}_k\big)_{ij} := \mathbb{E}\big[(\mathbf{q}_{k,i}-\overline{\mathbf{q}}_{k,i})(\mathbf{q}_{k,j}-\overline{\mathbf{q}}_{k,j})\big], \quad i,j=1,2,\dots,N,
\end{equation}
where $\mathbb{E}[\cdot]$ denotes the expectation operator and $\mathbf{q}_{k,i}$ denotes the $i$-th entry of the vector $\mathbf{q}_k$ (and similarly for the index $j$). Given an initial covariance matrix $\mathbf{P}_0$, the propagation of the covariance matrix from time $t_k$ to $t_{k+1}$ then follows from \eqref{eq:state prop}\,--\,\eqref{eq:covariance matrix def},
\begin{equation}\label{eq:cov prop}
	\mathbf{P}_{k+1} = \mathbf{M}_{k,k+1}\mathbf{P}_k\mathbf{M}_{k,k+1}\transpose,
\end{equation}
where superscript T the transpose. Equation \eqref{eq:cov prop} defines \emph{full-rank covariance propagation} and is the time propagation of the covariance matrix in the standard Kalman filter \citep[][Ch.~4 and references therein]{kalman1960new,simon2006optimal}. The variance $\bs{\sigma}^2_k$ at time $t_k$ is defined as the diagonal of the covariance matrix $\mathbf{P}_k$,
\begin{equation}\label{eq:variance def}
	\bs{\sigma}^2_k := \text{diag}(\mathbf{P}_k),
\end{equation}
and the correlation matrix $\mathbf{C}_k$ at time $t_k$ is the covariance matrix $\mathbf{P}_k$ normalized by the square root of the variances (standard deviations),
\begin{equation}\label{eq:correlation matrix def}
	(\mathbf{C}_k)_{ij} := \frac{(\mathbf{P}_{k})_{ij}}{\sqrt{\bs{\sigma}^2_{k,i}\bs{\sigma}^2_{k,j}}}, \quad i,j=1,2,\dots, N.
\end{equation}

We can compare the mean state and covariance propagation, \eqref{eq:mean state prop} and \eqref{eq:cov prop} respectively, with that of ensemble-based schemes, such as the EnKF. Ensemble-based data assimilation schemes begin with an ensemble of initial states $\mathbf{q}_0^m$ for $m=1,2,\dots,n_e$ sampled from a distribution with initial mean $\mathbf{q}_0$ and initial covariance $\mathbf{P}_0$ (such as a Gaussian distribution, for instance), and they evolve this ensemble forward in time according to the discrete state propagation defined in \eqref{eq:state prop}. At a time $t_k$, the exact mean state $\overline{\mathbf{q}}_k$ and exact covariance $\mathbf{P}_k$ are estimated by the ensemble mean $\overline{\mathbf{q}}_k^s$ and ensemble covariance $\mathbf{P}^s_k$ using the standard empirical estimators,
\begin{gather}
	\overline{\mathbf{q}}_k^s := \frac{1}{n_e}\sum_{m=1}^{n_e}\mathbf{q}_k^m,\label{eq:sample mean}\\
	\mathbf{P}_k^s := \frac{1}{n_e-1}\sum_{m=1}^{n_e}(\mathbf{q}_k^m-\overline{\mathbf{q}}_k^s)(\mathbf{q}_k^m-\overline{\mathbf{q}}_k^s)\transpose.\label{eq:sample cov}
\end{gather}
The ensemble variance $\bs{\sigma}^{2,s}_k$ can either be defined as the diagonal of the sample covariance matrix $\mathbf{P}_k^s$ or, equivalently, directly from the ensemble using the (unbiased) empirical estimator,
\begin{equation}\label{eq:sample variance}
	\bs{\sigma}^{2,s}_k := \frac{1}{n_e-1}\sum_{m=1}^{n_e}(\mathbf{q}_k^m-\overline{\mathbf{q}}_k^s)^2.
\end{equation}
The ensemble correlation matrix $\mathbf{C}_k^s$ can then be defined using \eqref{eq:correlation matrix def} by replacing the exact covariance and variance with the ensemble estimates in \eqref{eq:sample cov} and \eqref{eq:sample variance}, respectively.

Since the propagation of the individual ensemble members satisfies \eqref{eq:state prop}, it follows directly from \eqref{eq:state prop}--\eqref{eq:covariance matrix def} that the sample mean $\overline{\mathbf{q}}_k^s$ and sample covariance $\mathbf{P}_k^s$ satisfy the same propagation equations as the mean state $\overline{\mathbf{q}}_k$ and covariance $\mathbf{P}_k$, respectively, 
\begin{gather}
	\overline{\mathbf{q}}_{k+1}^s = \mathbf{M}_{k,k+1}\overline{\mathbf{q}}_k^s,\label{eq:sample mean prop} \\
	\mathbf{P}_{k+1}^s = \mathbf{M}_{k,k+1}\mathbf{P}_k^s\mathbf{M}_{k,k+1}\transpose. \label{eq:sample cov prop}
\end{gather}
Thus, the evolution of the mean state $\overline{\mathbf{q}}_k$ and covariance $\mathbf{P}_k$  are approximated through the evolution of the individual ensemble members. Since the number $n_e$ of ensemble members is typically much less than the dimension $N$ of the state, the sample covariance $\mathbf{P}_k^s$ is a low-rank approximation of the full-rank covariance $\mathbf{P}_k$. It can be shown that in the limit as the ensemble size $n_e$ approaches infinity, the approximated evolution of the mean state and covariance by the ensemble, \eqref{eq:sample mean prop} and \eqref{eq:sample cov prop} respectively, converge to the full rank propagation for the mean state, \eqref{eq:mean state prop}, and covariance \eqref{eq:cov prop} (\citealp[p.~231]{furrer2007estimation}; \citealp{butala2008asymptotic}; \citealp[p.~44]{evensen2009ensemble}). 

While full-rank covariance propagation \eqref{eq:cov prop} and ensemble-based approximations \eqref{eq:sample cov prop} both propagate the covariance, the means of propagation can be interpreted differently. Full-rank covariance propagation, as implemented in a Kalman filter, propagates the covariance matrix directly with its own, explicit evolution equation, in parallel with the evolution of the state. Ensemble-based methods, in contrast, do not evolve the covariance matrix explicitly per se, but rather do so indirectly through the propagation of the individual ensemble members.  From this perspective, one may expect that ensemble-based schemes, which do not evolve the covariance matrix directly, may not suffer from the same errors observed in full-rank covariance propagation. In the numerical experiments, we show that for advective dynamics that ensemble-based covariance propagation does in fact suffer from these same problems.

\section{Numerical Experiments}\label{sec:experiments}
The goal of the following numerical experiments is to illustrate that for a given discrete dynamical model, represented by the matrix $\mathbf{M}_{k,k+1}$, while the ensemble-based propagation of the mean state is accurate (up to the order of the numerical scheme), the ensemble-based covariance propagation can be remarkably inaccurate, far beyond that expected from sampling errors alone. We show that the ensemble-based covariance propagation faithfully captures the behavior produced by full-rank propagation, but neither approximates the exact covariance dynamics well. For these experiments, we consider states that satisfy versions of the advection equation for which the exact covariance dynamics are known. Thus, we can explicitly quantify the errors in both the ensemble-based and full-rank propagation schemes.

We start by defining the dynamical model for the state dynamics, then introduce the continuum framework from which we derive the corresponding continuum covariance dynamics and related quantities. We then consider two specific dynamical systems for the state, define the numerical discretizations, and describe the ensemble initialization. This is followed by the results of the experiments.

\subsection{Problem Setup}
Motivated by atmospheric data assimilation, we consider continuum states $q = q(x,t)$ for $x\in\mathbb{S}^1_1$, where $\mathbb{S}_1^1$ is the unit circle, that satisfy the following generalized advection equation,
\begin{equation}\label{eq:continuum state}
\begin{split}
	q_t + vq_x + bq = 0,\\
	q(x,t_0) = q_0(x),\,\,
	\end{split}
\end{equation}
where the velocity $v = v(x,t)$ and scalar $b = b(x,t)$ are deterministic. For our experiments, we will fix two choices of $b$ that define two different types of advective dynamics motivated by specific applications in atmospheric and chemical constituent data assimilation.

\subsubsection{Continuum Framework}\label{sec:continuum framework}
For states that satisfy the generalized advection equation \eqref{eq:continuum state}, the exact continuum covariance dynamics are also known. Therefore, we can leverage this information to evaluate and interpret the errors observed in full-rank and ensemble-based discrete covariance propagation. The full derivation and analysis of the continuum covariance dynamics are given in \cite{cohn1993dynamics} and \cite{gilpin2022continuum}, and we summarize the important concepts here.

Analogous to the discrete case, \eqref{eq:cov prop}, we can define the continuum covariance $P = P(x_1,x_2,t)$ for spatial variables $x_1,x_2 \in \mathbb{S}^1_1$ on the unit circle, associated with continuum states $q$,
\begin{equation}\label{eq:continuum covariance def}
	P(x_1,x_2,t) := \mathbb{E}\left[\{q(x_1,t)-\overline{q}(x_1,t)\}\{q(x_2,t)-\overline{q}(x_2,t)\}\right],
\end{equation}
where $\overline{q}(x_i,t) := \mathbb{E}[q(x_i,t)]$ denotes the mean state for $i=1,2$. Observing that the state dynamics are linear and coefficients $v$ and $b$ are deterministic, the continuum mean state $\overline{q}$ satisfies the same PDE as $q$, \eqref{eq:continuum state}. Following the definition of the covariance in \eqref{eq:continuum covariance def}, $P = P(x_1,x_2,t)$ satisfies its own PDE,
\begin{equation}\label{eq:continuum covariance pde}
\begin{split}
	P_t + v_1P_{x_1} + v_2P_{x_2} + (b_1+b_2)P = 0,\\
	P(x_1,x_2,t_0) = P_0(x_1,x_2), \quad\quad\,\,\,
\end{split}
\end{equation}
where the subscripts $1,2$ denote the quantities evaluated with respect to the spatial variables $x_1$ and $x_2$ \citep[see][for further details]{cohn1993dynamics,gilpin2022continuum}. We refer to this equation as the continuum covariance equation. We can relate \eqref{eq:continuum covariance pde} to the discrete, full-rank covariance propagation in \eqref{eq:cov prop} by interpreting  \eqref{eq:continuum covariance pde} as the evolution equation for the kernel $P=P(x_1,x_2,t)$ of the covariance operator $\bs{\mathcal{P}}_t$,
\begin{equation}
	(\bs{\mathcal{P}}_t f)(x_1) := \int_{\mathbb{S}^1_1}P(x_1,x_2,t)f(x_2)dx_2, \quad f \in L^2(\mathbb{S}_1^1).
\end{equation}
The covariance operator $\bs{\mathcal{P}}_t$ is defined on the space of square-integrable functions on the unit circle,\endnote{Specifically, $f \in L^2(\mathbb{S}_1^1)$ means $||f||_2 := \left(\int_{\mathbb{S}_1^1}|f(x)|^2dx\right)^{1/2} < \infty$.} $L^2(\mathbb{S}_1^1)$.
In addition, we define the fundamental solution operator $\bs{\mathcal{M}}_{t_k,t_{k+1}}$, which is an integral operator from $L^2(\mathbb{S}^1_1)$ to $L^2(\mathbb{S}^1_1)$ defined by \eqref{eq:continuum state} that propagates states from time $t_k$ to time $t_{k+1}$. With $\bs{\mathcal{M}}_{t_k,t_{k+1}}$ and its adjoint $\bs{\mathcal{M}}_{t_k,t_{k+1}}^*$ \citep[see][Sec.~2.1]{gilpin2022continuum}, the evolution of the covariance operator from time $t_k$ to $t_{k+1}$ is given by 
\begin{equation}\label{eq:continuum mpmt}
	\bs{\mathcal{P}}_{t_{k+1}} = \bs{\mathcal{M}}_{t_k,t_{k+1}}\bs{\mathcal{P}}_{t_k}\bs{\mathcal{M}}_{t_k,t_{k+1}}^*.
\end{equation}
Full-rank (discrete) covariance propagation defined in \eqref{eq:cov prop} is thus the discretization of \eqref{eq:continuum mpmt}.\endnote{Using the same fundamental solution operator, one can define the propagation of the continuum state $q$ from time $t_k$ to $t_{k+1}$ as
	$q(x,t_{k+1}) = \bs{\mathcal{M}}_{t_{k},t_{k+1}}q(x,t_k)$. Thus, discrete state propagation in \eqref{eq:state prop} is a discretization of this continuum equation.}

For nonzero initial correlation lengths (see \cite{gilpin2022continuum} and Sec.~\ref{sec:discussion} for further details), we can define the continuum variance,
\begin{equation}\label{eq:continuum variance def}
	 \sigma^2=\sigma^2(x,t) := \mathbb{E}\big[\{q(x,t)-\overline{q}(x,t)\}^2\big],
\end{equation}
which satisfies its own PDE,
\begin{equation}\label{eq:continuum variance}
\begin{split}
	\sigma^2_t + v\sigma^2_x + 2b\sigma^2 = 0,\\
	\sigma^2(x,t_0) = \sigma^2_0(x).\,\,\,
	\end{split}
\end{equation}
It then follows from the variance-correlation decomposition of the covariance,
\begin{equation}\label{eq:cov decomp}
	P(x_1,x_2,t) = \sigma(x_1,t)C(x_1,x_2,t)\sigma(x_2,t),
\end{equation}
that the continuum correlation $C = C(x_1,x_2,t)$ satisfies
\begin{equation}\label{eq:continuum correlation pde}
\begin{split}
	C_t + v_1C_{x_1} + v_2C_{x_2} = 0\\
	C(x_1,x_2,t_0) = C_0(x_1,x_2). 
	\end{split}
\end{equation}

From the correlations, we can define a correlation length scale $L(x,t)$,
\begin{equation}
L^2(x,t) := \big(-C_2(x,t)\big)^{-1}, \label{eq:correlation length def}
\end{equation}
where
\begin{equation}\label{eq:c2 definition}
 C_2(x,t) = \big[C_{x_1x_1}(x_1,x_2,t) - 2C_{x_1x_2}(x_1,x_2,t) + C_{x_2x_2}(x_1,x_2,t)\big]_{x_1=x_2=x}.
\end{equation}
This correlation length $L=L(x,t)$ satisfies its own PDE \citep[as derived in][]{cohn1993dynamics, gilpin2025inaccuracy},
\begin{equation}\label{eq:correlation length pde}
\begin{split}
L_t + vL_x - v_xL = 0,\\
L(x,t_0) = L_0(x).
\end{split}
\end{equation}
Observe that for a spatially-varying velocity field ($v_x \neq 0$) the correlation length field will evolve over space and time, even if the initial correlation length $L_0(x)$ is constant in space. In particular, $1/L$ satisfies the continuity equation, and therefore the average inverse correlation length is conserved in time.\endnote{I.e, the $L^1(\mathbb{S}_1^1)$ norm, $||1/L||_1 = \int_{\mathbb{S}_1^1}|1/L(x,t)|dx$, remains constant in time.} This implies that regions of convergence in the velocity field (where $v_x<0$) will cause correlation lengths to shrink, and will be balanced by regions of divergence in the velocity field (where $v_x > 0$) that will cause correlation lengths to grow. Thus, small correlation lengths are to be expected dynamically. We will see in the numerical experiments that these dynamics play a role during covariance propagation. 

For these numerical experiments, we let the velocity $v$ vary in space, but not in time,
\begin{equation}\label{eq:velocity}
	v(x) = \sin(x)+2.
\end{equation}
In this case, we have exact solutions to all previously stated continuum equations (for the state, covariance, variance, correlation, and correlation length) given in \citealp[Appen.~B of][]{gilpin2022continuum} and \citealp[Appen.~C of][]{gilpin2025inaccuracy}.

\subsubsection{Energy Conservation}\label{sec:energy conservation}
Motivated by chemical constituent data assimilation, we consider \eqref{eq:continuum state} with $b=\frac12v_x$. These dynamics arise when $q$ is the square root of a positive quantity, such as concentration or density. In particular, we rewrite these dynamics into conservation form,
\begin{equation}\label{eq:continuum state half}
\begin{split}
	q_t + \frac12(vq)_x + \frac12vq_x = 0,\\
	q(x,t_0) = q_0(x),\quad\,\,
\end{split}
\end{equation}
We refer to these dynamics as the  \emph{energy conservation case}, as the ``energy" defined as $\frac12||q(x,t)||_2^2$, where $||\cdot||_2$ defines the $L^2$ norm over $\mathbb{S}_1^1$, is conserved in time. The corresponding continuum covariance and variance equation satisfy conservation properties of their own, see Sec.~2 of \cite{gilpin2025inaccuracy} for more details.

We discretize \eqref{eq:continuum state half} using the Crank-Nicolson finite difference scheme, which applies a second-order centered-difference approximation of the spatial derivatives followed by the trapezoidal rule time-integration scheme \citep{crank1947practical}. The propagation of the discrete state $\mathbf{q}_k$ to time $t_{k+1}$ is defined as
\begin{equation}\label{eq:crank nicolson state}
\begin{split}
	\mathbf{q}_{k+1,i} + \frac{\Delta t}{8\Delta x}\left(\mathbf{v}_{i+1}\mathbf{q}_{k+1,i+1} - \mathbf{v}_{i-1}\mathbf{q}_{k+1,i-1}\right) + \frac{\Delta t}{8\Delta x}\mathbf{v}_i\left(\mathbf{q}_{k+1,i+1} - \mathbf{q}_{k+1,i-1}\right) \\= \mathbf{q}_{k,i} - \frac{\Delta t}{8\Delta x}\left(\mathbf{v}_{i+1}\mathbf{q}_{k,i+1} - \mathbf{v}_{i-1}\mathbf{q}_{k,i-1}\right) + \frac{\Delta t}{8\Delta x}\mathbf{v}_i\left(\mathbf{q}_{k,i+1} - \mathbf{q}_{k,i-1}\right),
	\end{split}
\end{equation}
where we assume periodic boundary conditions.
Crank-Nicolson is a second order (in space and time) accurate implicit scheme that is quadratically conservative when the velocity $v = v(x)$ does not depend on time, as in our numerical experiments \eqref{eq:velocity}. Therefore, discrete versions of the continuum conservation properties of the state, variance, and covariance, are preserved \citep[see Sec.~2 of][]{gilpin2025inaccuracy}. Since the Crank-Nicolson finite difference scheme is an implicit scheme, we do not explicitly compute the propagation matrix $\mathbf{M}_{k,k+1}$ and instead apply a linear solve to \eqref{eq:crank nicolson state} at each time step. Similarly, the corresponding full-rank covariance propagation, \eqref{eq:cov prop}, is not implemented directly since we do not form $\mathbf{M}_{k,k+1}$ explicitly. Instead, one time step corresponds to applying two linear solves, one application to compute $\tilde{\mathbf{P}}_{k+1} = \mathbf{M}_{k,k+1}\mathbf{P}_k$, which applies \eqref{eq:crank nicolson state} along the columns of $\mathbf{P}_{k}$, then $\mathbf{P}_{k+1} = \mathbf{M}_{k,k+1}\tilde{\mathbf{P}}_{k+1}\transpose$, which applies \eqref{eq:crank nicolson state} to the rows of $\tilde{\mathbf{P}}_{k+1}$.

The unit circle is discretized uniformly into $N=200$ grid points, $x_j = j\Delta x, \ \Delta x = 2\pi/N,\ j=0,1,2,\dots,N-1$. The time step $\Delta t$ is defined by the Courant number $\lambda$,
\begin{equation}\label{eq:cfl}
\lambda = \frac{\Delta t}{\Delta x}||v||,
\end{equation} 
where $||v||=3$ is the maximum value of the velocity field in \eqref{eq:velocity}.
While the Crank-Nicolson scheme is stable for any value of $\lambda$, we take $\lambda=1$ to define the time step $\Delta t$. We run these experiments to a final time $T_f = 3.98$ which corresponds to slightly after a full time period for these dynamics.

\subsubsection{Pure Advection}\label{sec:pure advection}
In addition to the energy conserving case, we consider the pure advection case where $b=0$, which has been studied in the chemical constituent data assimilation literature \citep{menard2000aassimilation,menard2021numerical,menard2000bassimilation,lyster2004lagrangian,pannekoucke2016parametric,pannekoucke2021methodology,pannekoucke2021anisotropic,sabathier2023boundary}. To remain consistent with those studies, we apply a first-order upwind forward-Euler finite difference scheme to define the discrete state propagation, 
\begin{equation}\label{eq:upwind fe}
	\mathbf{q}_{k+1,i} = \mathbf{q}_{k,i} - \frac{\Delta t}{\Delta x}\mathbf{v}_i\left(\mathbf{q}_{k,i}-\mathbf{q}_{k,i-1}\right).
\end{equation}
This is equivalent to a first-order semi-Lagrangian discretization scheme with linear interpolation. 
It is an explicit scheme, therefore we form the propagation matrix $\mathbf{M}_{k,k+1}$ explicitly to define the state propagation in \eqref{eq:state prop} and full-rank covariance propagation \eqref{eq:cov prop}. The spatial grid and time step $\Delta t$ are defined in the same way as for the energy conservation case. 

As noted in \cite{gilpin2025inaccuracy}, the upwind forward-Euler scheme has a transient phase in which numerical dissipation dominates before settling towards a roughly time-periodic behavior. Therefore, for the pure advection case we run experiments to $T_f=4.975$, slightly longer than the energy conservation case, so that solutions are no longer in the transient regime.

\subsubsection{Ensemble Initialization}
The initial ensemble $\{\mathbf{q}_0^m\}, \ m=1,2,\dots,n_e=4000$ is drawn from a multivariate Gaussian with constant initial mean state $\overline{\mathbf{q}}_0=4$ and initial covariance matrix $\mathbf{P}_0$. Initial covariances $\mathbf{P}_0$ are constructed using the Gaspari-Cohn compactly-supported approximation to a Gaussian \citep[][their (4.10)]{gaspari1999construction}. These correlation functions depend on a single cut-off length parameter $c$ that determines their region of compact support, and which is related to the correlation length $L(x,t)$ of \eqref{eq:correlation length def} by $L = \frac12\sqrt{0.3}c$. In our results, we use the initial value $c$ to distinguish between the covariances in the various experiments. We take the initial variance to be one in all experiments; experiments with a spatially-varying initial variance produced similar results to the unit variance case and are not shown. Initial ensemble members with negative values are discarded and replaced with new members that do not violate positivity. The Crank-Nicolson scheme \eqref{eq:crank nicolson state} and upwind forward Euler \eqref{eq:upwind fe} schemes are not positivity preserving. We do not enforce positivity during propagation and allow for the relatively few instances where ensemble members become negative during propagation.

\subsection{Results}\label{sec:results}
In this section, we present the results from the numerical experiments and highlight important features. Interpretation of these results is then given in Sec.~\ref{sec:discussion}.

\subsubsection{Mean State and Variances}\label{sec:mean and variance results}
\begin{figure}[!ht] 
   \centering
   \includegraphics[width=0.85\linewidth]{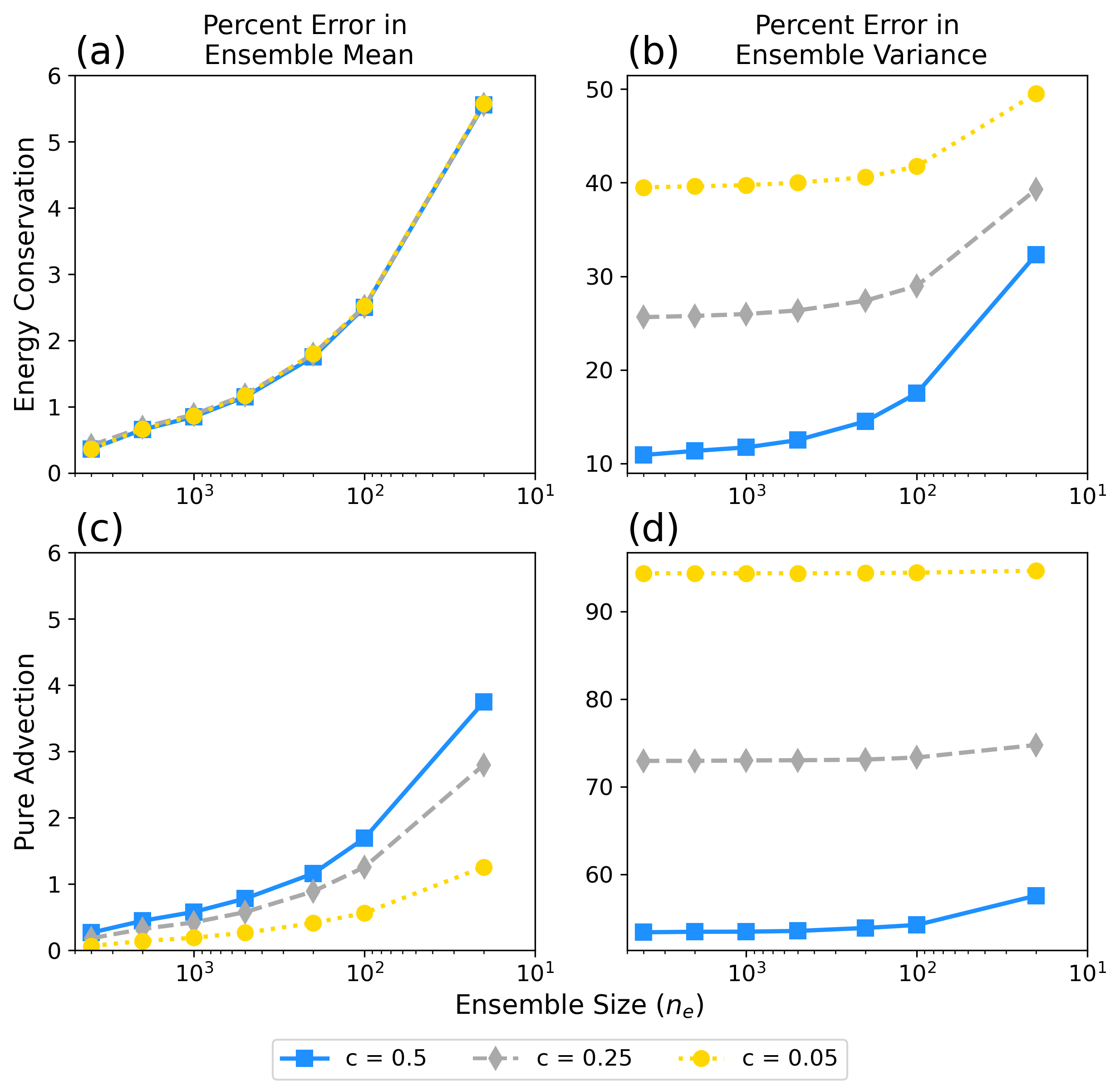} 
   \caption{\textit{Percentage error in the ensemble mean and ensemble variance as function of ensemble size and initial correlation length}. Panels (a) and (c) plot the normalized errors in the ensemble mean relative to the exact mean state as a function of ensemble size (horizontal axis) and initial correlation length (color and shape) for the energy conserving case (a) and pure advection case (c). Panels (b) and (d) plot the normalized error in the ensemble variance relative to the exact variance as a function of ensemble size (horizontal axis) and initial correlation length (color and shape) for the energy conserving (b) and pure advection (d) case. Errors are computed at the final time $T_f$ for the two respective numerical experiments, see text for details.} 
   \label{fig:figure1}
\end{figure}

Figure~\ref{fig:figure1} plots the average percent error in the $\ell^2$-norm\endnote{For an $N$-vector $\mathbf{x}$, the error in the $\ell^2$ norm with respect to $\tilde{\mathbf{x}}$ is defined as $||\mathbf{x}-\tilde{\mathbf{x}}||_{\ell^2} = (\sum_{i=1}^N(x_i-\tilde{x}_{i})^2)^{1/2}$, where $x_i$ and $\tilde{x}_{i}$ are the $i$-th components of $\mathbf{x}$ and $\tilde{\mathbf{x}}$, respectively.} in the ensemble mean (panels a and c) and ensemble variance (panels b and d) as a function of ensemble size and initial correlation length, denoted by the initial cut-off $c$ in the legend. Figures~\ref{fig:figure1}(a,b) correspond to the energy conservation case, and Figures~\ref{fig:figure1}(c,d) correspond to the pure advection case. In addition to the large ensemble $n_e=4000$, we calculate errors for ensemble sizes of $\{20,100,200,500,1000,2000\}$. For each fixed ensemble size, we randomly draw from the $n_e=4000$ ensemble of propagated states and compute the $\ell^2$ error in the mean state and variance estimates. We repeat this 1000 times and average the error over these repetitions, which is plotted in Figure~\ref{fig:figure1}.

\begin{figure}[htp] 
   \centering
   \includegraphics[width=0.8\linewidth]{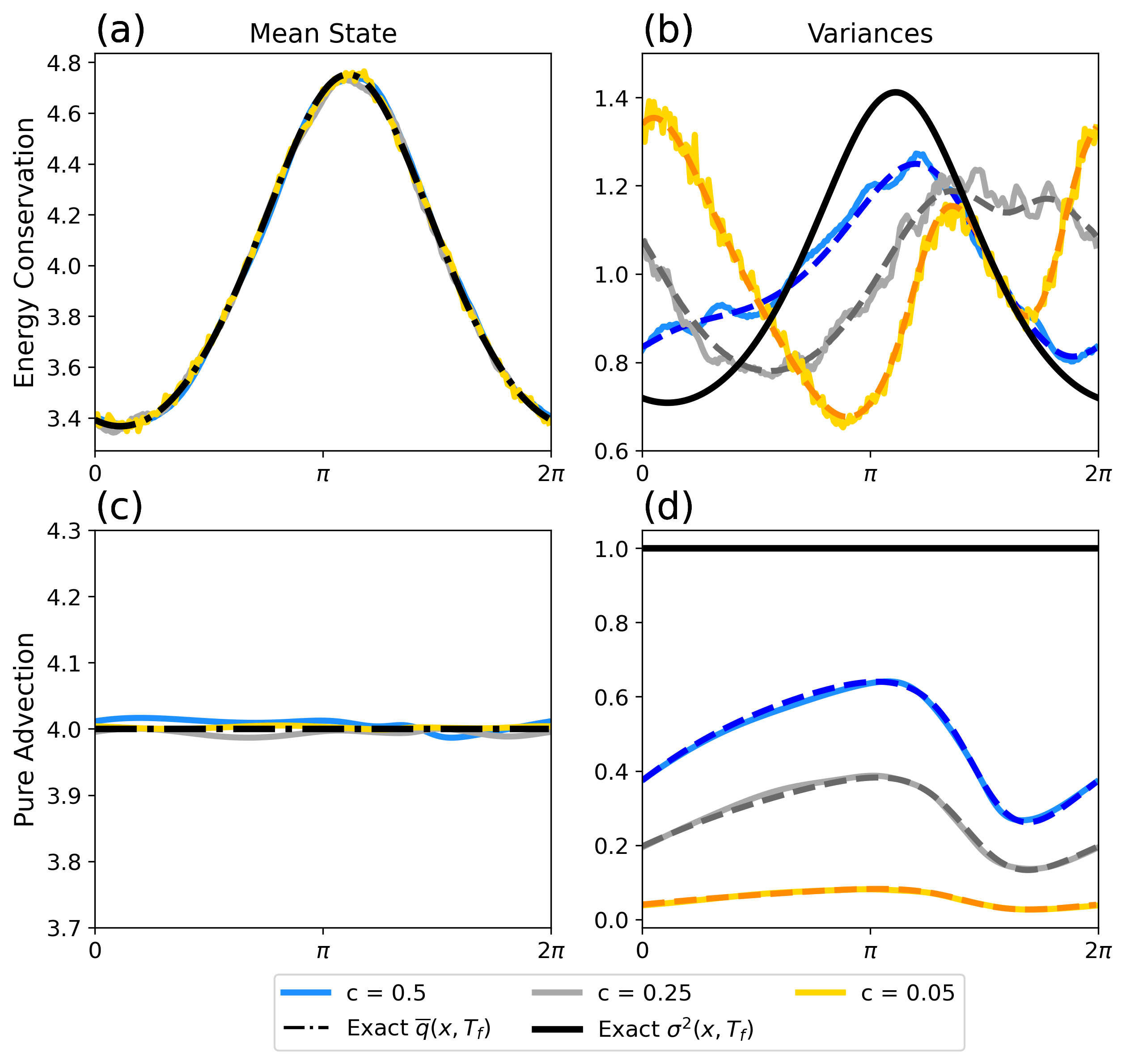} 
   \caption{\textit{Ensemble mean, ensemble variance, and covariance diagonals compared to the exact mean state and variance.} Panels (a) and (c) compare the exact mean state (denoted as Exact $\overline{q}(x,T_f)$ in the legend) with the ensemble mean for the different initial correlation lengths (corresponding to $c$ in the legend) for the energy conservation (a) and pure advection (c) cases. Panels (b) and (d) compare the exact variance (denoted as Exact $\sigma^2(x,T_f)$ in the legend) with the ensemble variance (solid) and variances approximated by full-rank covariance propagation (dashed) for different initial correlation lengths. The ensemble means and variances are calculated for ensemble size $n_e=4000$. Note the different vertical scales in each panel.} 
   \label{fig:figure2}
\end{figure}

The errors in the mean state estimated from the propagated ensemble, Figures~\ref{fig:figure1}(a,c),  are as expected: the errors are small (a few percent) and increase exponentially as the ensemble size decreases. In the energy conservation case (panel a), the errors do not depend on the initial correlation length, which is expected due to the conservation properties preserved by the Crank-Nicolson scheme. In the pure advection case (panel c), there is some dependence on the initial correlation length. This is likely due to errors in the numerical scheme as it approximates the exact mean state $q(x,t) = 4$, which remains constant for all time (see Figure~\ref{fig:figure2}c, for example). In general, the errors in the ensemble mean are small and within what is expected for these methods.

The behavior in the ensemble mean is in stark contrast to that of the ensemble variances, Figures~\ref{fig:figure1}(b,d). Across both cases, the errors in the ensemble variance are at least one order of magnitude larger than for the mean state, and they increase as the initial correlation length (proportional to $c$ in the legend) decreases. Perhaps surprisingly, the errors in the ensemble variances are relatively constant with respect to ensemble size, particularly for the pure advection case. We can see very clearly in Figures~\ref{fig:figure1}(b,d) that the errors in the ensemble variances depend primarily on correlation lengths, which is consistent with what is observed during full-rank covariance propagation \citep{gilpin2022continuum, gilpin2025inaccuracy}. This means that sampling errors are not the dominant source of error in the ensemble variances. Sampling errors only contribute significantly for the smallest ensemble (20 members), particularly in the pure advection case.

To better understand the contrast between the mean state and variance errors in Figure~\ref{fig:figure1}, Figure~\ref{fig:figure2} compares the exact mean state and variance with their ensemble approximations in the large ensemble case ($n_e=4000$), i.e., the mean and variance corresponding to the left-most points in Figure~\ref{fig:figure1}. For the variance comparison, we also plot the variances approximated by the corresponding full-rank covariance propagation (the diagonals from \eqref{eq:cov prop}). Figures~\ref{fig:figure2}(a,b) correspond to the energy conservation case, and Figures~\ref{fig:figure2}(c,d) the pure advection case. 

As expected from Figure~\ref{fig:figure1}, the ensemble means approximate the exact mean state rather well. The ensemble variances, however, are remarkably inaccurate. Figures~\ref{fig:figure2}(b,d) clearly show that the ensemble variances approximate the covariance diagonals extracted from full-rank covariance propagation very well, while neither are anywhere close to the exact variance dynamics. The variance approximations in the energy conservation case Figure~\ref{fig:figure2}(c) exhibit both loss and gain of variance. In fact, due to the conservation properties of the Crank-Nicolson scheme, the amount of variance lost is exactly equal to the amount of variance gained since the total variance is conserved in time. Thus, standard variance inflation, where the full ensemble variance is artificially inflated by a factor greater than one, may help in regions of variance loss, but not in regions of variance gain. In contrast, the pure advection case exhibits significant variance loss that becomes dramatically worse as correlation lengths decrease. Inflation schemes have been developed to address this loss of variance for this particular set of dynamics and numerical scheme \citep{menard2021numerical}, though it is not likely to generalize to higher-order schemes or to different dynamics, such as in the energy conservation case \citep[][Appen.~A]{gilpin2025inaccuracy}.

\subsubsection{Correlations}
\begin{figure}[!ht] 
   \centering
   \includegraphics[width=0.85\linewidth]{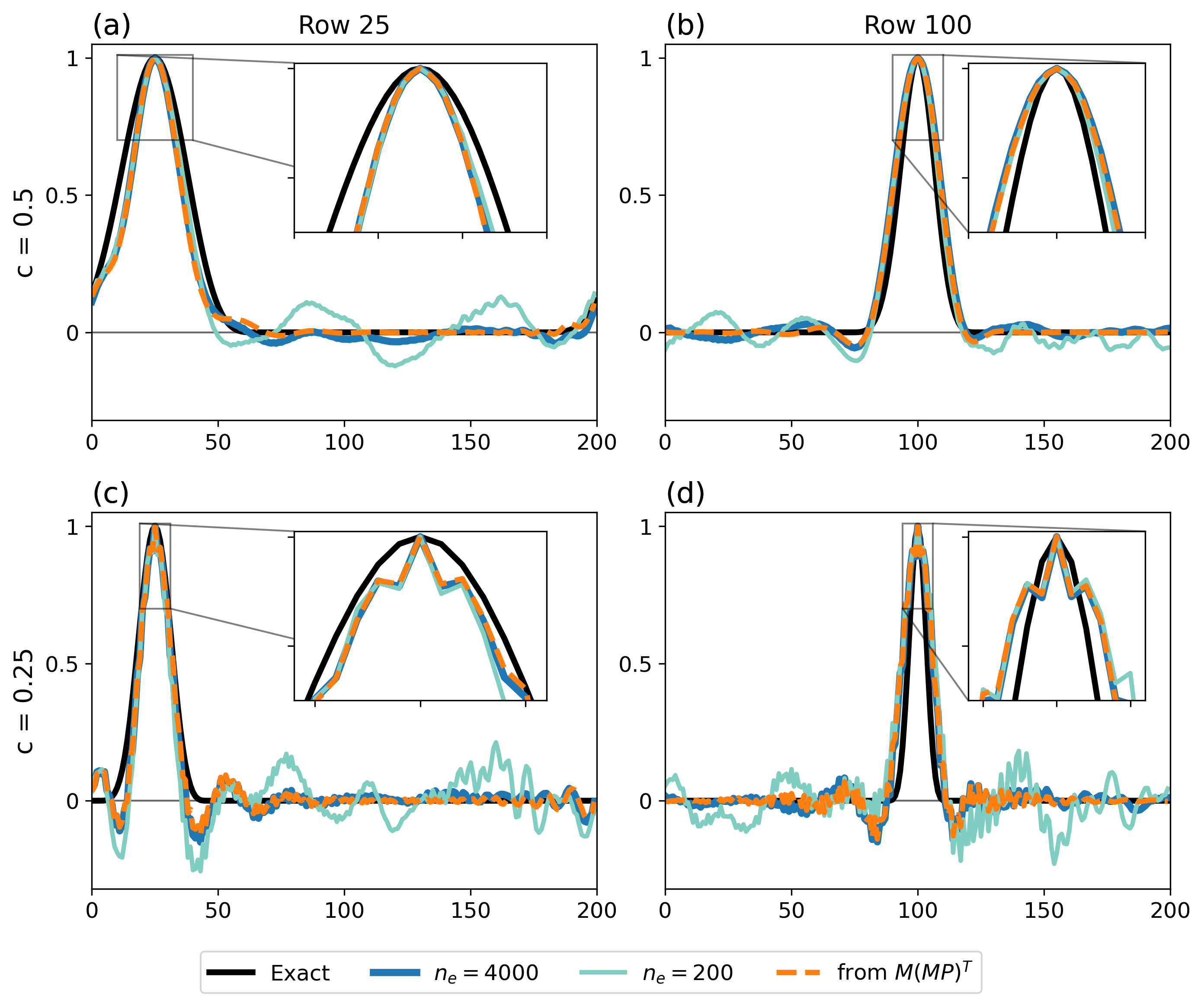} 
   \caption{\textit{One-dimensional correlations for the energy conservation case.} Each panel plots the exact correlation (solid black), correlation from full-rank covariance propagation (dashed orange), with ensemble correlations for $n_e=4000$ (solid blue) and $n_e=200$ (solid light teal) at the final time $T_f=3.98$ for the energy conservation dynamics. Panels (a) and (b) correspond to an initial covariance with $c=0.5$ and (c) and (d) with $c=0.25$. Row 25 of the correlation matrix is shown in (a) and (c); row 100 in (b) and (d). In each panel, the inset axis shows a zoomed in portion of the correlations near zero separation.} 
   \label{fig:figure3}
\end{figure}

\begin{figure}[ht] 
   \centering
   \includegraphics[width=0.85\linewidth]{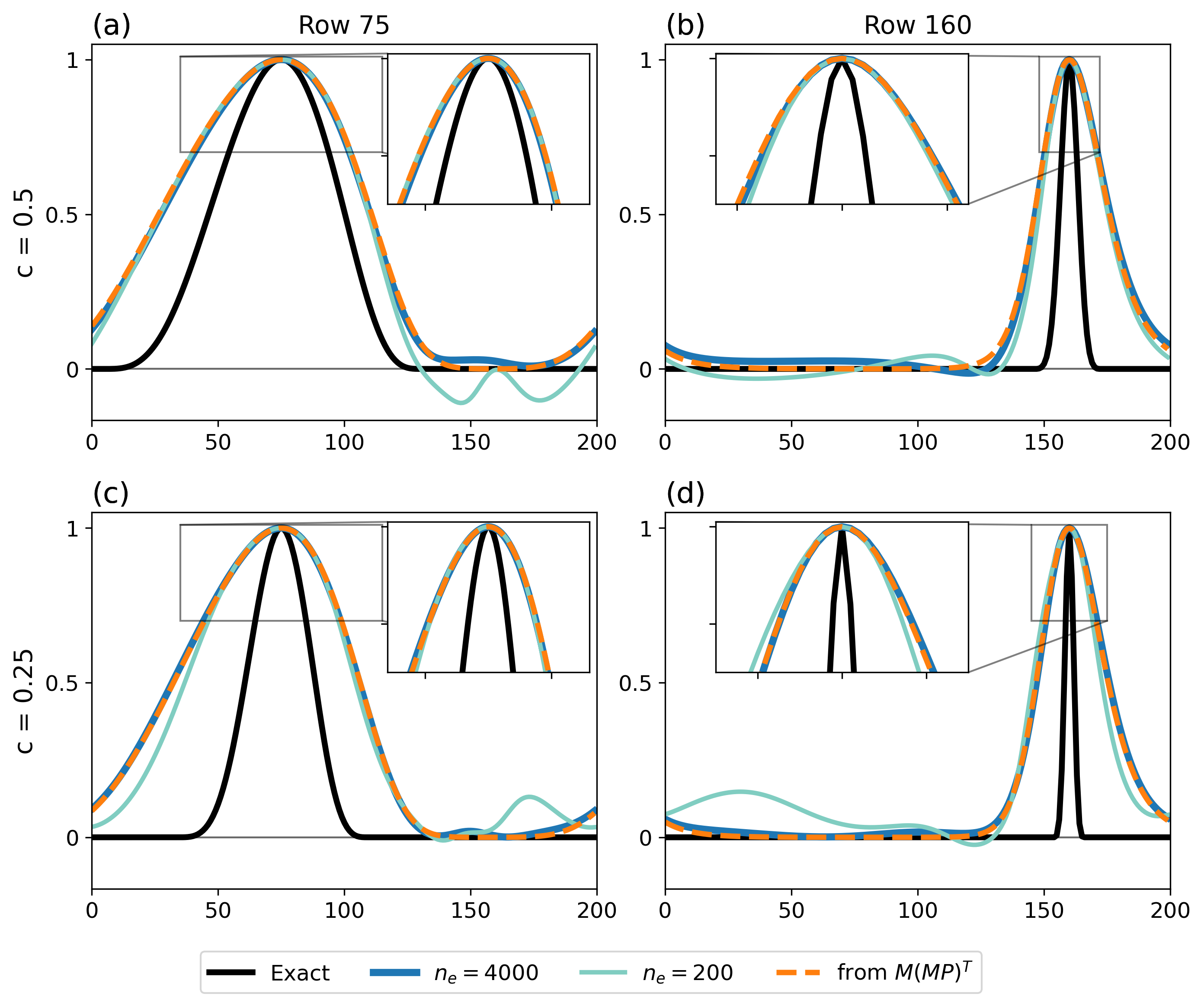} 
   \caption{\textit{One-dimensional correlations for the pure advection case.} Same as Figure~\ref{fig:figure3} for the pure advection case, corresponding to $T_f=4.975$ and for row 75 in panels (a) and (c) and row 160 in (b) and (d).} 
   \label{fig:figure4}
\end{figure}

In addition to the ensemble variances, we compare the exact correlations with correlations estimated from the ensemble and from full-rank covariance propagation. Figures~\ref{fig:figure3}--\ref{fig:figure4} show selected one-dimensional correlations approximated from the large ensemble ($n_e=4000$), a small ensemble ($n_e=200$), and correlations extracted from the corresponding full-rank covariance propagation. These are compared with the exact correlations for the energy conservation and pure advection cases, respectively. 

For both state dynamics, we see similar behavior in the ensemble correlations as seen in the ensemble variances. The ensemble correlations approximate correlations extracted from full-rank covariance propagation rather well, while neither can successfully approximate the exact correlations. Near zero separation, the ensemble correlations and those from full-rank covariance propagation are nearly identical. The errors in these approximations depend on the dynamics and location within the correlation matrix. In the energy conservation case, Figure~\ref{fig:figure3}, correlations are under-approximated in some regions, and over-approximated in other regions. The oscillatory behavior in the energy conservation case is due to numerical dispersion and is expected from the Crank-Nicolson scheme \citep[][p.~46]{beylkin1997adaptive}.

The correlations in the upwind case, Figure~\ref{fig:figure4}, are remarkably inaccurate near zero separation, even more so than in the energy conservation case in Figure~\ref{fig:figure3}. Initially, correlation lengths for the two cases $c=0.5$ and $c=0.25$ are small, but are fairly well resolved on the discrete grid. Over time, the approximated correlations increase dramatically and become almost indistinguishable between the $c=0.5$ and $c=0.25$ cases (compare down the columns of Figure~\ref{fig:figure4}). The increase in correlations is related to the variance loss observed in Figure~\ref{fig:figure1}(d) and Figure~\ref{fig:figure2}(d); this will be discussed further in Sec.~\ref{sec:discussion}.

Away from zero separation, spurious correlations in Figures~\ref{fig:figure3}--\ref{fig:figure4} begin to appear in the small ensemble case, which are not present in the large ensemble or full-rank approximations. Spurious, long range correlations are to be expected as ensemble sizes decrease, particularly when estimating small (in magnitude) correlations. Covariance localization, which tapers long-range correlations, is often applied in ensemble-based schemes to address this issue \citep{hamill2001distance,houtekamer2001sequential,hamill2009comments}. While covariance localization can mitigate spurious, long-range correlations due to sampling error, it will not address the errors in the correlations near zero separation.

\section{Discussion}\label{sec:discussion}
The numerical experiments in Sec.~\ref{sec:results} demonstrate the following:
\begin{enumerate}
	\item The propagation of the mean state from the ensemble is accurate, while the variances and correlations approximated by the same ensemble are remarkably inaccurate.
	\item Ensemble-based covariance propagation is a good approximation of full-rank covariance propagation, yet neither can approximate the exact covariance dynamics. 
\end{enumerate}
To understand why the covariances are so inaccurate in these numerical experiments requires taking a step back to understand the fundamental continuum aspects of covariance propagation.
  
On the surface, it may seem that the evolution of the covariance in ensemble-based covariance propagation and full-rank covariance propagation are fundamentally different. However, the two methods are linked together by 1) the same discrete dynamical model $\mathbf{M}_{k,k+1}$, and 2) the underlying continuum covariance dynamics, \eqref{eq:continuum covariance pde}. Since the ensemble-based and full-rank covariance propagation produce nearly identical results, as we have seen, we can interpret the experimental results in Sec.~\ref{sec:results} by studying the relationship between discrete covariance propagation, \eqref{eq:cov prop}, and the continuum covariance dynamics, \eqref{eq:continuum covariance pde}.

At an intuitive level, regions where correlation lengths become small correspond to regions with sharp gradients in the state, for example in regions of strong wind shear \citep{menard2000aassimilation, lyster2004lagrangian}. As a consequence one may expect the numerical schemes (i.e., that used to build $\mathbf{M}_{k,k+1}$), to struggle, as the accuracy of these methods begins to break down. 

Regions with sharp gradients, or equivalently short correlation lengths, are particularly problematic for advective dynamics, \eqref{eq:continuum state}, when the initial state is random and gives rise to the continuum covariance equation in \eqref{eq:continuum covariance pde}. As noted in \cite{gilpin2022continuum}, the hyperbolicity and the symmetry of the spatial derivatives in the continuum covariance equation \eqref{eq:continuum covariance pde} results in a peculiar behavior in the covariance dynamics that depends on the structure of the initial covariance, $P_0 = P_0(x_1,x_2)$. This peculiar behavior is that the continuum covariance equation has two distinct types of solutions. First, when $P_0$ has nonzero initial correlation lengths, the covariance solution $P=P(x_1,x_2,t)$ has nonzero correlation lengths for all time (i.e., remains spatially correlated), and the dynamics along $x_1=x_2$ satisfies the variance equation, \eqref{eq:continuum variance}. Second, if the initial covariance $P_0$ is spatially uncorrelated (i.e., strictly diagonal), the covariance $P$ remains spatially uncorrelated (diagonal) for all time, and the dynamics along $x_1=x_2$ satisfies instead what is called the continuous spectrum equation,
\begin{equation}\label{eq:continuous spectrum equation}
\begin{split}
P^c_t + vP^c_x + (2b-v_x)P^c = 0,\\
P^c(x,t_0) = P^c_0(x),\quad\quad\,\,
\end{split}
\end{equation}
where we refer to $P^c = P^c(x,t)$ as the continuous spectrum solution \citep[Sec.~2.3 of][]{gilpin2022continuum}.
Equivalently stated, white noise remains white. These two solutions to the continuum covariance equation  \eqref{eq:continuum covariance pde} are distinct when $v_x \neq 0$. Formally, this means there exists a discontinuity in the continuum covariance equation \eqref{eq:continuum covariance pde} in the limit as correlation lengths tend to zero.

While the initial covariance $P_0$ may not be truly uncorrelated in practice, the continuous spectrum solution makes its presence known during discrete covariance propagation, \eqref{eq:cov prop}, when correlation lengths become small. This behavior manifests in error terms which arise in the continuum covariance dynamics that \eqref{eq:cov prop} is approximating. The analysis of the continuum dynamics approximated by \eqref{eq:cov prop} is broken down into the behavior along the covariance diagonal presented in \cite{gilpin2025inaccuracy} and the behavior in the correlations presented here in Appen.~\ref{sec:appendix correlation}. The error terms in the approximated variance and correlation dynamics depend explicitly on powers of the ratio of the grid spacing $\Delta x$ and correlation length $L(x,t)$, the same correlation length which satisfies \eqref{eq:correlation length pde} (see (45) and (47) of \cite{gilpin2025inaccuracy} for the variances and \eqref{eq:cd correlation dynamics} and \eqref{eq:up correlation dynamics} in Appen~\ref{sec:appendix correlation} for the correlations). In fact, the form of the error terms in the approximated variance and correlation dynamics are nearly the same and in some cases mirror each other in their behavior. These error terms are different from the typical discretization errors expected for a numerical scheme (i.e., local truncation error \cite[pp. 104--105]{leveque1992numerical}) and are a direct consequence of \eqref{eq:cov prop} approximating across the diagonal $x_1=x_2$. Specifically, the inherent structure of \eqref{eq:cov prop} approximates across the discontinuity in the continuum covariance dynamics that distinguishes between the two solutions \citep[see][for further discussion]{gilpin2025inaccuracy}.

These error terms depending on powers of $\Delta x/L(x,t)$ that appear in the variance and correlation dynamics approximated by \eqref{eq:cov prop} explain the inaccurate behavior observed in the numerical experiments in Sec.~\ref{sec:results}. As correlation lengths shrink according to \eqref{eq:correlation length pde}, for a fixed grid resolution, these error terms can become large enough to completely alter the approximated dynamics. In the energy conservation case, the leading error terms are of the order $\Delta x^2/L^2$ (\citealp[see][their (45) for the variances]{gilpin2025inaccuracy}; \eqref{eq:cd correlation dynamics} in Appen~\ref{sec:appendix correlation} for the correlations). As correlation lengths become small, these terms produce the loss and gain of variance seen in Figure~\ref{fig:figure2}(b) and corresponding errors in the correlations of Figure~\ref{fig:figure3}. In the pure advection case, the first-order upwind spatial discretization yields a leading order error term of order $\Delta x/L^2(x,t)$ (\citealp[see][their (47) for the variances]{gilpin2025inaccuracy}; \eqref{eq:up correlation dynamics} in Appen~\ref{sec:appendix correlation} for the correlations). In the approximated variance dynamics, this error term can be interpreted as a zeroth order (with respect to derivative) dissipative term in the variance approximations causing significant variance loss. A corresponding error term of opposite sign occurs in the approximated correlation dynamics, mirroring the effect on the variances, namely as variance is lost, this loss is compensated by an increase in correlations. This relationship between variance loss on increased correlations has been observed \citep[][]{menard2000aassimilation}, and is now clarified by identifying these explicit error terms. 

While these error terms depend on the grid resolution $\Delta x$, increasing the spatial resolution (decreasing $\Delta x$) is not likely to ameliorate the errors observed during discrete covariance propagation. The initial correlation length scales in these numerical experiments are well resolved, and both shrink and grow over space and time according to \eqref{eq:correlation length pde}. In general, regions of wind shear, i.e., sharp gradients in the velocity field, will continually produce smaller correlation lengths, cascading into smaller scales \citep[pp. 2327--2328]{lyster2004lagrangian}. As correlation lengths continue to shrink, finer grid resolution (smaller $\Delta x$) will be necessary to ensure these error terms $\Delta x/L(x,t)$ produced by \eqref{eq:cov prop} remain small. Taking into account the computational cost of decreasing $\Delta x$ and the dynamics of the correlation length scales, increasing spatial resolution may not lead to a simple fix for the variance and correlations.

To see that ensemble-based covariance propagation is almost identical to full-rank covariance propagation is both surprising and not surprising. One would expect that as the number of ensemble members $n_e$ becomes large, the ensemble covariances should approximate the covariances  propagated at full rank by \eqref{eq:cov prop}. This is verified in our numerical experiments. What is perhaps surprising is that the propagation of the mean state by the ensemble is accurate while the ensemble covariance is not, all while using the same discrete dynamical model $\mathbf{M}_{k,k+1}$ defined by the state dynamics. These simple numerical experiments, together with previous analysis, bring to light a potentially fundamental problem for data assimilation with states that satisfy advective dynamics. The numerical experiments demonstrate that one cannot simply assume the propagation of the covariance, whether directly by \eqref{eq:cov prop}, or indirectly by an ensemble, will be taken care of by appropriately propagating the state. This phenomenon is due to the hyperbolicity of advective dynamics, and is likely not entirely the case when hyperbolicity is lost, for instance in the presence of a diffusive term \citep[][pp.~893--894]{gilpin2022continuum}.  While the continuum state and covariance dynamics are related, their behavior as PDEs are quite different and should therefore be treated as such. Covariance propagation schemes, as they currently stand, are not suited for covariance propagation for advective dynamics, and new methods should be explored with this behavior in mind.

\section{Conclusions}\label{sec:conclusions}
Through a series of simple numerical experiments, we have demonstrated that accurate mean state propagation by an ensemble does not imply accurate covariance propagation when correlation lengths become small. For states that satisfy the linear advection equation and its variants, the exact continuum covariance dynamics are known, and we have compared the mean state and covariance estimated from a propagated ensemble of states with the exact mean state and covariance. From these experiments, we observe the following:
\begin{itemize}
	\item While the mean state estimates are fairly accurate, the errors in the ensemble variances are an order of magnitude larger than for the state, they are fairly constant as the ensemble size decreases, and they increase as initial correlation lengths decrease.
	\item The ensemble covariances converge to the covariances approximated by full-rank covariance propagation, and neither approximate the exact covariance well.
	\item The primary source of error in ensemble-based covariance propagation is due to the discrete covariance propagation itself, much more so than typical sampling errors due to finite ensemble size.
\end{itemize}

In particular, these numerical experiments demonstrate that ensemble-based covariance propagation produces covariances that are nearly identical to those produced during full-rank covariance propagation. This implies that while each propagation scheme may seem different, the underlying problem is what connects the two, namely the discrete dynamical model $\mathbf{M}_{k,k+1}$ and the continuum covariance dynamics \eqref{eq:continuum covariance pde}. The analysis of the errors in full-rank covariance propagation in \cite{gilpin2025inaccuracy} for the variances and for the correlations here in Appen.~\ref{sec:appendix correlation} clearly identify the error terms that produce the inaccurate discrete covariance dynamics. These terms involve powers of the ratio of the grid resolution $\Delta x$ and the correlation length $L(x,t)$, which become large as correlation lengths become small. The resulting errors in both the full-rank and ensemble-based covariances are severe, and for the ensemble-based covariances in particular, cannot be rectified by the usual inflation and localization schemes.

While these numerical experiments are rather simple, the results are profound and enlightening. Ensemble-based data assimilation is widely practiced in the data assimilation community. While the discrete mean state propagation is understood to be accurate, at least within the limitations of grid resolution, the corresponding discrete covariance propagation is also assumed to be accurate and can therefore be overlooked. These simple numerical experiments illustrate that accurate discrete covariance propagation is not necessarily the case. The root of the problem lies in the relationship between the continuum covariance dynamics and the discrete covariance propagation defined from the discrete state. This work indicates that covariance propagation for advective dynamics can be problematic across data assimilation schemes that evolve covariances, independent of the propagation scheme. Thus, a shift in mindset regarding covariance propagation in data assimilation may be necessary. Such a shift is already starting to take form in methods such as local covariance evolution introduced in \cite{cohn1993dynamics} and \cite{gilpin2023new} and the parametric Kalman filter \citep{pannekoucke2016parametric,pannekoucke2021sympkf,perrot2022toward}.

\section*{Data Accessibility}
The code used to run the numerical experiments and generate the figures in this article are available from the corresponding author upon request.

\section*{Acknowledgements}
The author would like to acknowledge Dr. Stephen E. Cohn and Dr. Tomoko Matsuo for their time and expertise that helped to strengthen and clarify the contributions of this work.

\section*{Funding Information}
SG is supported in part by the Data Driven Discovery RTG at the University of Arizona under NSF Grant No. DMS-1937229. 

\section*{Competing Interests}
The author has no competing interest to declare.

\section*{Authors' Contributions}
The author contributed to the conceptualization, analysis, investigation, software development, visualization, writing, and funding acquisition. 

\appendix
\section{Analysis of the Correlation Dynamics}\label{sec:appendix correlation}
Figures~\ref{fig:figure3}--\ref{fig:figure4} show examples of the correlations approximated by the ensemble and from the full-rank covariance propagation for the two numerical experiments. These experiments demonstrate that in addition to the variances, the correlations are poor approximations of the exact correlations. In this appendix, we provide an error analysis of the discretized correlations approximated by \eqref{eq:cov prop} in order to identify the terms that cause the errors in these approximations. This analysis follows the same procedure in \cite{gilpin2025inaccuracy}, which derives the continuum dynamics being approximated by \eqref{eq:cov prop} along the covariance diagonal. Here, we do the same for the correlations. 

Although the results presented in Figures~\ref{fig:figure3}--\ref{fig:figure4} correspond to a full discretization, (i.e., the discretization of the continuum state dynamics in both space and time), it is sufficient to study the semi-discretization, namely discretization only in space, while leaving time a continuous variable. For this analysis, we will assume that the correlation $C=C(x_1,x_2,t)$ is four times continuously differentiable in both spatial arguments, and the velocity $v=v(x,t)$ is at least four times continuously differentiable in its spatial argument.

Recall that correlation $C=C(x_1,x_2,t)$ for $x_1,x_2 \in \mathbb{S}_1^1$ satisfies \eqref{eq:continuum correlation pde}. Observe from \eqref{eq:continuum correlation pde} that the correlation dynamics are independent of $b$ in the state dynamics, and for $C_0 \geq0$ (as is the case in our numerical experiments), $C(x_1,x_2,t) \geq 0$ for all $t\geq t_0$.

We discretize the unit circle $\mathbb{S}_1^1$ such that $x_i = i\Delta x$ for $i=1,2,\dots,N$, $\Delta x = 2\pi/N$, and assume periodicity. Quantities evaluated on the spatial grid are denoted by
\begin{equation}
	q(x_i,t) := q_i, \quad v(x_i,t) := v_i,\quad P(x_i,x_j,t) := P_{ij},\quad C(x_i,x_j,t) := C_{ij}.
\end{equation}

Given a spatial discretization for the continuum state dynamics, we can derive the corresponding semi-discretization for the correlations, i.e., the time evolution of the elements $C_{ij}$. To obtain this semi-discretization, we start with the semi-discretization for the covariance $P_{ij}$, then apply the following decomposition,
\begin{equation}\label{eq:squareroot covariance}
	P_{ij} := \sqrt{P_{ii}}C_{ij}\sqrt{P_{jj}}, \quad i,j=1,2,\dots N,
\end{equation}
where $P_{ii}$ and $P_{jj}$ denote the $i$-th and $j$-th diagonal entries of the covariance. 

We consider two cases corresponding to the schemes used in the numerical experiments in Sec.~\ref{sec:results} and the analysis of the covariance diagonals in \cite{gilpin2025inaccuracy}. First, we consider a first-order upwind spatial discretization applied to the generalized advection equation \eqref{eq:continuum state} for general $b$, noting that $b=0$ for pure advection is a special case of this analysis. Second, we consider a second-order centered difference discretization applied to the energy conservation dynamics in \eqref{eq:continuum state half}. To ease the derivations and notation, we present the analysis in the reverse order of the numerical experiments, starting with the upwind differencing followed by centered differencing.

\subsection{First-Order Upwind Differencing}
The pure advection case in Sec.~\ref{sec:pure advection} is a special case of the generalized advection equation, \eqref{eq:continuum state}. Therefore, we will consider the generalized advection equation for this analysis (noting that taking $b=0$ corresponds to pure advection). Consider \eqref{eq:continuum state} and apply a first-order upwind discretization to the spatial derivative $q_x$ to obtain the semi-discretization for the state,
\begin{equation}
	\ddt q_i = \frac{1}{\Delta x}v_i(q_{i-1}-q_i) - b_iq_i.
\end{equation}
To define the corresponding semi-discretization for the covariance, we first define the covariance $P_{ij}$ in terms of the error $\varepsilon$,
\begin{equation}
\varepsilon_i := q_i - \mathbb{E}\left[{q}_i\right], \quad P_{ij} := \mathbb{E}\left[\varepsilon_i\varepsilon_j\right].
\end{equation}
The semi-discretization for the covariance then follows,
\begin{equation}\label{eq:upwind pij}
	\ddt P_{ij} = \frac{1}{\Delta x}v_i\left[P_{i-1,j} - P_{ij}\right] +  \frac{1}{\Delta x}v_j\left[P_{i,j-1} - P_{ij}\right] - (b_i+b_j)P_{ij}.
\end{equation}

To derive the semi-discretization for the correlation, replace $P_{ij}$ in \eqref{eq:upwind pij} with \eqref{eq:squareroot covariance}, expand the time derivative on the left-hand side, and isolate $\frac{d}{dt}C_{ij}$,
\begin{align}\label{eq:upwind cij 1}
	\ddt C_{ij} &= -\frac12C_{ij}\left[\frac{1}{P_{ii}}\ddt P_{ii} + \frac{1}{P_{jj}}\ddt P_{jj}\right] \nonumber\\
	&+ \frac{v_i}{\Delta x}\left[\sqrt{\frac{P_{i-1,i-1}}{P_{ii}}}C_{i-1,j} - C_{ij}\right]
	+ \frac{v_j}{\Delta x}\left[\sqrt{\frac{P_{j-1,j-1}}{P_{jj}}}C_{i,j-1} - C_{ij}\right]
	-(b_1+b_j)C_{ij}.
\end{align}
From \eqref{eq:upwind pij}, we have expressions for $\frac{1}{P_{ii}}\ddt P_{ii}$ and $\frac{1}{P_{jj}}\ddt P_{jj}$, which we can substitute into \eqref{eq:upwind cij 1}. After a bit of simplification, we arrive at the semi-discretization for the correlation for first-order upwind differencing,
\begin{align}\label{eq:upwind cij 2}
	\ddt C_{ij} &= \frac{v_i}{\Delta x}\sqrt{\frac{P_{i-1,i-1}}{P_{ii}}}\bigg[ C_{i-1,j} - \frac12C_{ij}\left(C_{i-1,i}+C_{i,i-1}\right)\bigg]\nonumber\\
	& \ +  \frac{v_j}{\Delta x}\sqrt{\frac{P_{j-1,j-1}}{P_{jj}}}\bigg[ C_{i,j-1} - \frac12C_{ij}\left(C_{j-1,j}+C_{j,j-1}\right)\bigg].
\end{align}
Note that for $C_{ii}=1, \ i=1,2,\dots,N$ initially, the cancelation in the bracketed terms implies that $C_{ii}=1,  \ i=1,2,\dots,N$ for all time.
Observe also that while this semi-discretization is defined for the full correlation matrix, the second set of terms inside the square brackets are local terms that average across the covariance matrix diagonal. Even though these averaging terms are local, they impact the whole correlation matrix.

The goal now is to derive the continuum dynamics approximated by \eqref{eq:upwind cij 2}. To do so, we first introduce a new quantity through a change of variables $x = (x_1+x_2)\slash2$, $\xi = (x_1-x_2)\slash2$,
\begin{equation}
	\tilde{C} = \tilde{C}(x,\xi,t) := C\big(x_1(x,\xi), x_2(x,\xi),t\big),
\end{equation}
then expand $\tilde{C}(x,\xi,t)$ about $\xi=0$, which yields,
\begin{equation}
	\tilde{C}(x,\xi,t) = C_0(x,t) + \frac{\xi^2}{2}C_2(x,t) + \mathcal{O}(\xi^4) =  1+ \frac{ \xi^2}{2}C_2(x,t) + \mathcal{O}(\xi^4),
\end{equation}
where we assume that $C_0(x,t) = C(x,x,t) =1$. Note that $C_2$ is the same quantity defined in \eqref{eq:c2 definition} \citep[see e.g.,][pp.~3136--3137]{cohn1993dynamics}. We can also define $P_0(x,t) = P(x,x,t)$ in the same way as $C_0$ above. We use this change of variables to expand the pairs of local averaging terms inside the square brackets of \eqref{eq:upwind cij 2}. These terms are second-order centered approximations of the correlation on the half grid, specifically,
\begin{equation}
	C_{i-1,i}+C_{i,i-1} = \tilde{C}\big(x_{i-1/2},-\frac{\Delta x}{2},t\big) + \tilde{C}\big(x_{i-1/2},\frac{\Delta x}{2},t\big) = 2 + \frac{\Delta x^2}{4}C_2(x_{i-1/2},t) + \mathcal{O}(\Delta x^4),\label{eq:average minus}
\end{equation}
and similarly for index $j$. We also use this change of variables and $P_0$ to simplify the square root quantities. Considering the first term in \eqref{eq:upwind cij 2}, expanding the numerator about $x_i$ yields
\begin{align}
	\sqrt{\frac{P_{i-1,i-1}}{P_{i,i}}} &= 1-\frac12\Delta x \big(\log P_0(x_i,t)\big)_x + \mathcal{O}(\Delta x^2).\label{eq:square root minus}
\end{align}

Making these substitutions and with a bit of rearranging,
\begin{align}
	\ddt C_{ij} &= \frac{v_i}{\Delta x}\big[ C_{i-1,j} - C_{ij}\big] + \frac{v_j}{\Delta x}\big[ C_{i,j-1} - C_{ij}\big] \nonumber\\
	&-\frac{v_i}{2}\big(\log P_0(x_i,t)\big)_x\big[ C_{i-1,j} - C_{ij}\big] -\frac{v_j}{2}\big(\log P_0(x_j,t)\big)_x\big[ C_{i,j-1} - C_{ij}\big]\nonumber \\
	& \ + C_{ij}\frac{\Delta x}{8}\left[\frac{v_i}{L^2(x_{i-1/2},t)} + \frac{v_j}{L^2(x_{j-1/2},t)}\right]\label{eq:upwind2 line3}\\
	& - \frac{\Delta x^2}{16}C_{ij}\left[\frac{v_i\big(\log P_0(x_i,t)\big)_x}{L^2(x_{i-1/2},t)} + \frac{v_j\big(\log P_0(x_j,t)\big)_x}{L^2(x_{j-1/2},t)}\right] 
	+ \mathcal{O}(\Delta x^3)\nonumber
\end{align}
We have replaced $C_2(x_{i-1/2},t)$ and $C_2(x_{j-1/2},t)$ with the correlation length using \eqref{eq:correlation length def}.
By expanding these terms back to $x_i$ and $x_j$, we have the continuum correlation dynamics approximated by first-order upwind differencing,
\begin{align}\label{eq:up correlation dynamics}
	\ddt C(x_i,x_j,t) &= -v(x_i,t)C_{x_i}(x_i,x_j,t)-v(x_j,t)C_{x_j}(x_i,x_j,t)\nonumber\\
	&\ + \frac{\Delta x}{2}\left[v(x_i,t)C_{x_ix_i}(x_i,x_j,t) + v(x_j,t)C_{x_jx_j}(x_i,x_j,t)\right]\nonumber\\
	&\ + \frac{\Delta x}{2}\left[v(x_i,t)\big(\log P_0(x_i,t)\big)_xC_{x_i}(x_i,x_j,t) + v(x_j,t)\big(\log P_0(x_j,t)\big)_xC_{x_j}(x_i,x_j,t)\right]\nonumber\\
	&\ + \frac{\Delta x}{8}C(x_i,x_j,t)\left[\frac{v(x_i,t)}{L^2(x_i,t)} + \frac{v(x_j,t)}{L^2(x_j,t)}\right] \\
	&\ - \frac{\Delta x^2}{16}C(x_i,x_j,t)\bigg[v(x_i,t)\left(\frac{1}{L^2(x_i,t)}\right)_x + \frac{v(x_i,t)\big(\log P_0(x_i,t)\big)_x}{L^2(x_i,t)}\nonumber\\ 
	& \quad + v(x_j,t)\left(\frac{1}{L^2(x_j,t)}\right)_x + \frac{v(x_j,t)\big(\log P_0(x_j,t)\big)_x}{L^2(x_j,t)}\bigg] + G_u(x_i,x_j,t) + H_u(x_i,x_j,t).\nonumber
\end{align}
Here, $H_u(x_i,x_j,t)$ is $\mathcal{O}(\Delta x^2)$, involving higher-order spatial derivatives of the correlation in $x_i$ and $x_j$. The term $G_u(x_i,x_j,t)$ is $\mathcal{O}(\Delta x^3)$ and contains higher-order error terms in powers of $\Delta x/L$.

The first two terms on the right-hand side of \eqref{eq:up correlation dynamics} yield the correct correlation dynamics, \eqref{eq:continuum correlation pde}. The first and second set of bracketed terms on the right-hand side of order $\Delta x$ are the typical, first-order error terms expected from first-order upwind differencing. The third and fourth set of bracketed terms on the right-hand side of \eqref{eq:up correlation dynamics} arise from the local averaging terms in the square brackets of \eqref{eq:upwind cij 2}. These error terms depend on powers of the ratio $\Delta x/L(x,t)$ and can become significant as correlation lengths become small for fixed $\Delta x$. 

The third set of bracketed terms on the right-hand side of \eqref{eq:up correlation dynamics} can be considered as a first-order error term that causes correlations to increase when $C$ is positive, and to decrease when $C$ is negative. In our numerical experiments, the initial correlation $C_0 \geq0$, and therefore by the correlation PDE \eqref{eq:continuum correlation pde}, $C(x_1,x_2,t)\geq 0$ for all time. Since $L>0$ and $v>0$ by assumption, this term in our experiments is positive. Thus, correlations will grow over time, and faster for small correlation lengths. This error term is of the opposite sign and half the magnitude of an identical error term in the corresponding approximated dynamics along the covariance diagonal, Eq.~(47) of \cite{gilpin2025inaccuracy}. The mirroring of this error term in the correlation and covariance diagonal yields two conclusions. First, this balancing of decreasing variance and increasing correlation is consistent with the conservation properties of the continuum correlation dynamics. Second, the primary source of error in the correlations is caused by \eqref{eq:cov prop} averaging across the covariance diagonal.

The fourth set of bracketed terms on the right-hand side of \eqref{eq:up correlation dynamics} is also a consequence of the local averaging terms in \eqref{eq:upwind2 line3}, though its impact on the correlation dynamics is generally less than the first-order error term that precedes it. This term is important when comparing these approximated dynamics with those from the centered differencing case discussed in the next section.

\subsection{Second-Order Centered Differencing}
The Crank-Nicolson finite difference scheme in Sec.~\ref{sec:energy conservation} applies a second-order centered difference spatial discretization to \eqref{eq:continuum state half}. Therefore, we proceed here in the same way, applying second-order centered differencing to the spatial derivatives $(vq)_x$ and $q_x$ in \eqref{eq:continuum state half} to define the semi-discretization for the state,
\begin{equation}
	\ddt q_i = \frac{1}{4\Delta x}\left(v_{i-1}q_{i-1} - v_{i+1}q_{i+1}\right) + \frac{1}{4\Delta x}v_i\left(q_{i-1}-q_{i+1}\right)
\end{equation}
Following the same procedure outlined in the first-order upwind case, we arrive at the semi-discretization for the covariance,
\begin{align}
	\ddt P_{ij} &= \frac{1}{4\Delta x}\left(v_{i-1}P_{i-1,j} - v_{i+1}P_{i+1,j}\right)+ \frac{1}{4\Delta x}v_i\left(P_{i-1,j} - P_{i+1,j}\right)\nonumber \\
	& \, + \frac{1}{4\Delta x}\left(v_{j-1}P_{i,j-1} - v_{j+1}P_{i,j+1}\right) + \frac{1}{4\Delta x}v_j\left(P_{i,j-1} - P_{i,j+1}\right),
\end{align}
and the corresponding semi-discretization for the correlation,
\begin{align}
	\ddt C_{ij} &= \frac{1}{2\Delta x}\sqrt{\frac{P_{i-1,i-1}}{P_{ii}}}\left(\frac{v_{i-1}+v_i}{2}\right)\left[C_{i-1,j} - \frac12C_{ij}(C_{i-1,i}+C_{i,i-1})\right]\nonumber \\ 
	& \ - \frac{1}{2\Delta x}\sqrt{\frac{P_{i+1,i+1}}{P_{ii}}}\left(\frac{v_{i+1}+v_i}{2}\right)\left[C_{i+1,j} - \frac12C_{ij}(C_{i+1,i}+C_{i,i+1})\right]\nonumber \\
	& \ + \frac{1}{2\Delta x}\sqrt{\frac{P_{j-1,j-1}}{P_{jj}}}\left(\frac{v_{j-1}+v_j}{2}\right)\left[C_{i,j-1} - \frac12C_{ij}(C_{j-1,j}+C_{j,j-1})\right]\label{eq:cd corr sd} \\ 
	& \ - \frac{1}{2\Delta x}\sqrt{\frac{P_{j+1,j+1}}{P_{jj}}}\left(\frac{v_{j+1}+v_j}{2}\right)\left[C_{i,j+1} - \frac12C_{ij}(C_{j+1,j}+C_{j,j+1})\right]\nonumber
\end{align}
Observe that the semi-discretization in the centered difference case is quite similar to the upwind case. The second term in each of the square brackets is the same type of local averaging term seen in the upwind case that produces correlation quantities along the half grid. In addition, by pairing the first and second line terms and the third and fourth line terms together, these reflect flux-like spatial derivatives also  on the half grid, noting that the velocity terms are grouped in such a way to make this averaging on the half grid more clear. 

We derive the continuum dynamics approximated by \eqref{eq:cd corr sd} using the same procedure applied to the upwind case with a few additional terms. The pairs of local averaging terms within each of the four square brackets are expanded following \eqref{eq:average minus} and 
\begin{equation}\label{eq:average plus}
	C_{i+1,i}+C_{i,i+1} = \tilde{C}\big(x_{i+1/2},\frac{\Delta x}{2},t\big) + \tilde{C}\big(x_{i+1/2},-\frac{\Delta x}{2},t\big) = 2 + \frac{\Delta x^2}{4}C_2(x_{i+1/2},t) + \mathcal{O}(\Delta x^4).
\end{equation}
There are two sets of square root terms that need to be replaced, where one set can be replaced with \eqref{eq:square root minus}, and the other with the following,
\begin{align}
	\sqrt{\frac{P_{i+1,i+1}}{P_{ii}}} 
	 &= 1+\frac12\Delta x \big(\log P_0(x_i,t)\big)_x + \mathcal{O}(\Delta x^2)\label{eq:square root plus},
\end{align}
as in \eqref{eq:square root minus}. In addition, the velocity terms represent (second-order) centered approximations of velocities along the half grid, which we see in the following expansions,
\begin{align}
	\frac{v_{i-1}+v_i}{2} &= v(x_{i-1/2},t) + \frac{\Delta x^2}{8} v_{xx}(x_{i-1/2},t) + \mathcal{O}(\Delta x^4),\\
	\frac{v_{i+1}+v_i}{2} &= v(x_{i+1/2},t) + \frac{\Delta x^2}{8} v_{xx}(x_{i+1/2},t) + \mathcal{O}(\Delta x^4),
\end{align}
and similarly for index $j$. Making use of this expansion for the velocity and replacing the square root quantities in \eqref{eq:cd corr sd}, we arrive at
\begin{align}
\ddt C_{ij} &= \frac{1}{4\Delta x}\left\{ \bigg[v_{i-1}C_{i-1,j} - v_{i+1}C_{i+1,j}\big] + C_{ij}(v_{i+1}-v_{i-1}) + v_i\big[C_{i-1,j}-C_{i+1,j}]\right\}\nonumber \\
& \ +  \frac{1}{4\Delta x}\left\{ \bigg[v_{j-1}C_{i,j-1} - v_{j+1}C_{i,j+1}\big] + C_{ij}(v_{j+1}-v_{j-1}) + v_j\big[C_{i,j-1}-C_{i,j+1}]\right\}\nonumber \\
&-\frac18\big(\log P_0(x_i,t)\big)_x\left[(v_{i-1}+v_i)(C_{i-1,j}-C_{ij}) + (v_{i+1}+v_i)(C_{i+1,j}-C_{ij})\right]\label{eq:cd corr2}\\ 
&-\frac18\big(\log P_0(x_j,t)\big)_x\left[(v_{j-1}+v_j)(C_{i,j-1}-C_{ij}) + (v_{j+1}+v_j)(C_{i,j+1}-C_{ij})\right]\nonumber\\ 
& +C_{ij}\frac{\Delta x^2}{16}\bigg\{[v(x_i,t)C_2(x_i,t)]_x +[v(x_j,t)C_2(x_j,t)]_x\nonumber \\
&\hspace{1in} + v_i\big(\log P_0(x_i,t)\big)_xC_2(x_i,t) + v_j\big(\log P_0(x_j,t)\big)_xC_2(x_j,t) \bigg\} + \mathcal{O}(\Delta x^3)\nonumber
\end{align}
In the last step, we expand terms back to $x_i$ and $x_j$. We then obtain the continuum dynamics being approximated by the centered difference scheme,
\begin{align}\label{eq:cd correlation dynamics}
	\ddt C(x_i,x_j,t) &= -v(x_i,t)C_{x_i}(x_i,x_j,t) -v(x_j,t)C_{x_j}(x_i,x_j,t) \nonumber  \\
	& \ - \frac{\Delta x^2}{12}\bigg[3v_x(x_i,t)C_{x_ix_i}(x_i,x_j,t) +3v_{xx}(x_i,t)C_{x_i}(x_i,x_j,t) + 2v(x_i,t)C_{x_ix_ix_i}(x_i,x_j,t) \nonumber \\
	 & \ \quad +3v_x(x_j,t)C_{x_jx_j}(x_i,x_j,t) + 3v_{xx}(x_j,t)C_{x_j}(x_i,x_j,t) + 2v(x_j,t)C_{x_jx_jx_j}(x_i,x_j,t)\bigg] \nonumber \\
	 & \ -\frac{\Delta x^2}{4} \big(\log P_0(x_i,t)\big)_x\bigg[v(x_i,t)C_{x_ix_i}(x_i,x_j,t) + v_x(x_i,t)C_{x_i}(x_i,x_j,t)\bigg] \\
	 & \ -\frac{\Delta x^2}{4} \big(\log P_0(x_j,t)\big)_x\bigg[v(x_j,t)C_{x_jx_j}(x_i,x_j,t) + v_x(x_j,t)C_{x_j}(x_i,x_j,t)\bigg] \nonumber \\
	& \ - \frac{\Delta x^2}{16}C(x_i,x_j,t)\bigg[\left(\frac{v(x_i,t)}{L^2(x_i,t)}\right)_x + \frac{v(x_i,t)\big(\log P_0(x_i,t)\big)_x}{L^2(x_i,t)} \nonumber \\ 
	& \ + \left(\frac{v(x_j,t)}{L^2(x_j,t)}\right)_x + \frac{v(x_j,t)\big(\log P_0(x_j,t)\big)_x}{L^2(x_j,t)}\bigg] + G_c(x_i,x_j,t) + H_c(x_i,x_j,t).\nonumber
\end{align}
The term $H_c(x_i,x_j,t)$ is $\mathcal{O}(\Delta x^4)$ involving higher-order spatial derivatives of the correlation, velocity, and $P_0$ in the $x_i$ and $x_j$ directions. The term $G_c(x_i,x_j,t)$ is $\mathcal{O}(\Delta x^4)$ and contains higher-order terms involving powers of $\Delta x/L$. 

As observed in the upwind case, the first and second terms on the right-hand side of \eqref{eq:cd correlation dynamics} recover the correct correlation dynamics of \eqref{eq:continuum correlation pde}. 
The first, second, and third sets of terms in the square brackets are the typical, second-order error terms expected for the centered difference scheme.
The fourth set of bracketed terms on the right-hand side of \eqref{eq:cd correlation dynamics} is a direct consequence of the local averaging terms in \eqref{eq:cd corr sd}, and like the upwind case, depends explicitly on the ratio of the grid resolution $\Delta x$ and correlation length $L$. Thus, as correlation lengths become small, these error terms become large and corrupt the correlation propagation by \eqref{eq:cov prop}. Unlike the upwind case, where correlations grow over time, the sign of the fourth bracketed term on the right-hand side of  \eqref{eq:cd correlation dynamics} can change. This is reflected in Figure~\ref{fig:figure3}, where in the correlations from the ensemble and \eqref{eq:cov prop} are under-approximations in some regions and over-approximations in others. 

\subsection{Concluding Remarks}
The correlation dynamics being approximated by \eqref{eq:cov prop} for the upwind and centered difference schemes in \eqref{eq:up correlation dynamics} and \eqref{eq:cd correlation dynamics}, respectively, clarify the errors observed in the numerical experiments, Figures~\ref{fig:figure3} and \ref{fig:figure4}. The errors in the correlations are caused by the local averaging terms in \eqref{eq:upwind cij 2} and \eqref{eq:cd corr sd}, which produce terms in the approximated correlation dynamics that depend on the ratio of the grid resolution $\Delta x$ and correlation length $L$. As correlation lengths become small for fixed $\Delta x$, these error terms become large. In the upwind case, the term of order $\Delta x/L^2$ causes an increase in correlations for positive correlations, which mirrors the corresponding loss of variance along the diagonal. In the centered-difference case, the term of order $\Delta x^2/L^2$ results in over- or under-approximated correlations depending on the sign of the error terms. 

The equations for the approximated correlation dynamics presented here are consistent with those derived for the covariance diagonal in \cite{gilpin2025inaccuracy}. Together, they provide a more complete picture of the inaccuracy of discrete covariance propagation for advective systems. The underlying problem is the disconnect between the inherent structure of \eqref{eq:cov prop} (i.e., the local averaging terms) and the peculiar behavior of the continuum covariance dynamics in \eqref{eq:continuum covariance pde} (its two types of solutions). 
\theendnotes

\bibliography{../references}

\begin{thebibliography}{}

\bibitem[Beylkin and Keiser, 1997]{beylkin1997adaptive}
Beylkin, G. and Keiser, J. (1997).
\newblock An adaptive pseudo-wavelet approach for solving nonlinear partial
  differential equations.
\newblock In {\em Multiscale Wavelet Methods for Partial Differential
  Equations}, volume~6 of {\em Wavelet Analysis and Applications}. Academic
  Press.

\bibitem[Butala et~al., 2008]{butala2008asymptotic}
Butala, M.~D., Yun, J., Chen, Y., Frazin, R.~A., and Kamalabadi, F. (2008).
\newblock Asymptotic convergence of the ensemble {K}alman filter.
\newblock In {\em 2008 15th IEEE International Conference on Image Processing},
  pages 825--828. IEEE.

\bibitem[Cohn, 1993]{cohn1993dynamics}
Cohn, S. (1993).
\newblock Dynamics of short-term univariate forecast error covariances.
\newblock {\em Mon. Weather Rev.}, pages 3123--3150.

\bibitem[Crank and Nicolson, 1947]{crank1947practical}
Crank, J. and Nicolson, P. (1947).
\newblock A practical method for numerical evaluation of solutions of partial
  differential equations of the heat-conduction type.
\newblock {\em Math. Proc. Cambridge}, 43:50--67.

\bibitem[Daley, 1991]{daley1991atmospheric}
Daley, R. (1991).
\newblock {\em Atmospheric Data Analysis}.
\newblock Cambridge University Press.

\bibitem[Dee, 1995]{dee1995line}
Dee, D.~P. (1995).
\newblock On-line estimation of error covariance parameters for atmospheric
  data assimilation.
\newblock {\em Mon. Weather Rev.}, 123(4):1128--1145.

\bibitem[Evensen, 1994]{evensen1994sequential}
Evensen, G. (1994).
\newblock {Sequential data assimilation with a nonlinear quasi-geostrophic
  model using Monte Carlo methods to forecast error statistics}.
\newblock {\em J. Geophys. Res.}, 99:10143--10162.

\bibitem[Evensen, 2009]{evensen2009ensemble}
Evensen, G. (2009).
\newblock {\em {Data Assimilation: The Ensemble Kalman Filter}}.
\newblock Springer, 2nd edition.

\bibitem[Furrer and Bengtsson, 2007]{furrer2007estimation}
Furrer, R. and Bengtsson, T. (2007).
\newblock Estimation of high-dimensional prior and posterior covariance
  matrices in {K}alman filter variants.
\newblock {\em Journal of Multivariate Analysis}, 98(2):227--255.

\bibitem[Gaspari and Cohn, 1999]{gaspari1999construction}
Gaspari, G. and Cohn, S.~E. (1999).
\newblock Construction of correlation functions in two and three dimensions.
\newblock {\em {Q. J. Roy. Meteor. Soc.}}, 125:723--757.

\bibitem[Gilpin, 2023]{gilpin2023new}
Gilpin, S. (2023).
\newblock {\em {A new perspective on covariance propagation for data
  assimilation applications}}.
\newblock PhD thesis, University of Colorado Boulder.

\bibitem[Gilpin et~al., 2022]{gilpin2022continuum}
Gilpin, S., Matsuo, T., and Cohn, S.~E. (2022).
\newblock Continuum covariance propagation for understanding variance loss in
  advective systems.
\newblock {\em SIAM/ASA Journal on Uncertainty Quantification}, 10(3):886--914.

\bibitem[Gilpin et~al., 2025]{gilpin2025inaccuracy}
Gilpin, S., Matsuo, T., and Cohn, S.~E. (2025).
\newblock Inaccuracy of the variance evolution associated with discrete
  covariance propagation.
\newblock {\em {Q. J. Roy. Meteor. Soc.}}, pages 1--26.

\bibitem[Hamill et~al., 2009]{hamill2009comments}
Hamill, T.~M., Whitaker, J.~S., Anderson, J.~L., and Snyder, C. (2009).
\newblock Comments on ``{S}igma-point {K}alman filter data assimilation methods
  for strongly nonlinear systems''.
\newblock {\em Journal of the Atmospheric Sciences}, 66(11):3498 -- 3500.

\bibitem[Hamill et~al., 2001]{hamill2001distance}
Hamill, T.~M., Whitaker, J.~S., and Snyder, C. (2001).
\newblock Distance-dependent filtering of background error covariance estimates
  in an ensemble {K}alman filter.
\newblock {\em Mon. Weather Rev.}, 129(11):2776--2790.

\bibitem[Houtekamer and Mitchell, 1998]{houtekamer1998data}
Houtekamer, P.~L. and Mitchell, H.~L. (1998).
\newblock {Data assimilation using an ensemble Kalman filter technique}.
\newblock {\em Mon. Weather Rev.}, 126:769--811.

\bibitem[Houtekamer and Mitchell, 2001]{houtekamer2001sequential}
Houtekamer, P.~L. and Mitchell, H.~L. (2001).
\newblock {A sequential ensemble Kalman filter for atmospheric data
  assimilation}.
\newblock {\em Mon. Weather Rev.}, 129:123--137.

\bibitem[Jazwinski, 1970]{jazwinski1970stochastic}
Jazwinski, A.~H. (1970).
\newblock {\em Stochastic Processes and Filter Theory}.
\newblock Academic Press.

\bibitem[Kalman, 1960]{kalman1960new}
Kalman, R. (1960).
\newblock A new approach to linear filtering and prediction problems.
\newblock {\em Journal of Basic Engineering}, 82:35--45.

\bibitem[LeVeque, 1992]{leveque1992numerical}
LeVeque, R. (1992).
\newblock {\em Numerical Methods for Conservation Laws}.
\newblock Birkhauser Verlag.

\bibitem[Lorenc, 2003]{lorenc2003modelling}
Lorenc, A.~C. (2003).
\newblock Modelling of error covariances by 4d-var data assimilation.
\newblock {\em {Q. J. Roy. Meteor. Soc.}}, 129(595):3167--3182.

\bibitem[Lyster et~al., 2004]{lyster2004lagrangian}
Lyster, P.~M., Cohn, S.~E., Zhang, B., Chang, L.-P., M\'{e}nard, R., Olson, K.,
  and Renka, R. (2004).
\newblock {A Lagrangian trajectory filter for constituent data assimilation}.
\newblock {\em {Q. J. Roy. Meteor. Soc.}}, 130:2315--2334.

\bibitem[M\'{e}nard and Chang, 2000]{menard2000bassimilation}
M\'{e}nard, R. and Chang, L.~P. (2000).
\newblock {Assimilation of stratospheric chemical tracer observations using a
  Kalman Filter. Part II: $\chi^2$-validated results and analysis of variance
  and correlation dynamics}.
\newblock {\em Mon. Weather Rev.}, 128:2672--2686.

\bibitem[M\'{e}nard et~al., 2000]{menard2000aassimilation}
M\'{e}nard, R., Cohn, S.~E., Chang, L.~P., and Lyster, P.~M. (2000).
\newblock {Assimilation of stratospheric chemical tracer observations using a
  Kalman filter. Part I: Formulation}.
\newblock {\em Mon. Weather Rev.}, 128:2654--2671.

\bibitem[M\'{e}nard et~al., 2021]{menard2021numerical}
M\'{e}nard, R., Skachko, S., and Pannekoucke, O. (2021).
\newblock Numerical discretization causing error variance loss and the need for
  inflation.
\newblock {\em {Q. J. Roy. Meteor. Soc.}}, 147(740):3498--3520.

\bibitem[Pannekoucke, 2021]{pannekoucke2021anisotropic}
Pannekoucke, O. (2021).
\newblock An anisotropic formulation of the parametric kalman filter
  assimilation.
\newblock {\em Tellus A}, 73(1):1--27.

\bibitem[Pannekoucke and Arbogast, 2021]{pannekoucke2021sympkf}
Pannekoucke, O. and Arbogast, P. (2021).
\newblock Sympkf (v1. 0): a symbolic and computational toolbox for the design
  of parametric kalman filter dynamics.
\newblock {\em Geoscientific Model Development}, 14(10):5957--5976.

\bibitem[Pannekoucke et~al., 2021]{pannekoucke2021methodology}
Pannekoucke, O., M\'{e}nard, R., Aabaribaoune, M.~E., and Plu, M. (2021).
\newblock {A methodology to obtain model-error covariances due to the
  discretization scheme from the parametric Kalman filter perspective}.
\newblock {\em Nonlin. Processes Geophys.}, 28:1--22.

\bibitem[Pannekoucke et~al., 2016]{pannekoucke2016parametric}
Pannekoucke, O., Ricci, S., Barthelemy, S., M\'{e}nard, R., and Thual, O.
  (2016).
\newblock {Parametric Kalman Filter for chemical transport model}.
\newblock {\em Tellus A}, 68:31547.

\bibitem[Perrot et~al., 2022]{perrot2022toward}
Perrot, A., Pannekoucke, O., and Guidard, V. (2022).
\newblock Toward a multivariate formulation of the pkf assimilation:
  application to a simplified chemical transport model.
\newblock {\em EGUsphere}, 2022:1--43.

\bibitem[Sabathier et~al., 2023]{sabathier2023boundary}
Sabathier, M., Pannekoucke, O., Maget, V., and Dahmen, N. (2023).
\newblock Boundary conditions for the parametric kalman filter forecast.
\newblock {\em Journal of Advances in Modeling Earth Systems},
  15(10):e2022MS003462.

\bibitem[Simon, 2006]{simon2006optimal}
Simon, D. (2006).
\newblock {\em Optimal State Estimation: Kalman, $H_\infty$, and Nonlinear
  Approaches}.
\newblock Wiley.

\bibitem[Whitaker and Hamill, 2012]{whitaker2012evaluating}
Whitaker, J. and Hamill, T.~M. (2012).
\newblock {Evaluating methods to account for system errors in ensemble data
  assimilation}.
\newblock {\em Mon. Weather Rev.}, pages 3078--3090.

\end{thebibliography}
\end{document}